\documentclass[11pt]{article}

\usepackage[a4paper,top=2.5cm,bottom=2.5cm,left=2.5cm,right=2.5cm]{geometry}

\usepackage{amssymb}
\usepackage{amsmath}
\usepackage{amsthm}
\usepackage[english]{babel}
\usepackage[latin1]{inputenc}
\usepackage{latexsym}
\usepackage{mathrsfs}
\usepackage{bbm}
\usepackage{graphicx}
\usepackage[hidelinks]{hyperref}
\usepackage{accents}
\usepackage{cases}
\usepackage{url}
\usepackage{color}

\usepackage{chngcntr}
\counterwithout{equation}{section}

\newtheorem{Theorem}{Theorem}[section]
\newtheorem{Proposition}[Theorem]{Proposition}
\newtheorem{Definition}[Theorem]{Definition}
\newtheorem{Lemma}[Theorem]{Lemma}

\theoremstyle{remark} 
\newtheorem{remark}[Theorem]{Remark}

\newcommand{\nn}{\mathbb{N}}
\newcommand{\rr}{\mathbb{R}}

\newcommand{\eee}{\mathbb{E}}

\newcommand{\aaa}{\mathcal{A}}

\newcommand{\eps}{\varepsilon}

\newcommand{\h}{\alpha}
 
\newcommand{\R}{\mathbb{R}}

\newcommand{\E}{\mathcal{E}}
\newcommand{\F}{\mathcal{F}}

\DeclareMathOperator{\esssup}{esssup}
\def\esssup_#1{\underset{#1}{\mathrm{ess\,sup\, }}}
\def\essinf_#1{\underset{#1}{\mathrm{ess\,inf\, }}}
\def\argmax_#1{\underset{#1}{\mathrm{arg\,max\, }}}
\def\argmin_#1{\underset{#1}{\mathrm{arg\,min\, }}}

\def \ep{\hbox{ }\hfill$\Box$}
\def \Sum{\displaystyle\sum}

\def \Int{\displaystyle\int}

\def \Inf{\displaystyle\inf}

\def\b1{\bf 1}

\def \I{\mathbb{I}}

\def \R{\mathbb{R}}
\def \M{\mathbb{M}}

\def \E{\mathbb{E}}
\def \F{\mathbb{F}}

\def \P{\mathbb{P}}

\def \S{\mathbb{S}}

\def \Ac{{\cal A}}

\def \Cc{{\cal C}}
\def \Dc{{\cal D}}

\def \Fc{{\cal F}}
\def \Gc{{\cal G}}

\def \Sc{{\cal S}}

\def \Wc{{\cal W}}

\def \ep{\hbox{ }\hfill$\Box$}
 
\def\reff#1{{\rm(\ref{#1})}}
\def\beqs{\begin{eqnarray*}}
\def\enqs{\end{eqnarray*}}
\def\beq{\begin{eqnarray}}
\def\enq{\end{eqnarray}}

\def \trans{^{\scriptscriptstyle{\intercal}}}

\begin{document}

\title{A Weak Martingale Approach to Linear-Quadratic McKean-Vlasov Stochastic Control Problems
}

\author{
Matteo BASEI
\footnote{Industrial Engineering and Operations Research Department (IEOR), University of California, Berkeley, \sf basei at berkeley.edu}
\qquad\quad
Huy\^en PHAM
\footnote{Corresponding author. Laboratoire de Probabilités, Statistique et Modélisation (LPSM), Université Paris Diderot and CREST-ENSAE, \sf pham at lpsm.paris}
}

\maketitle

\renewcommand{\abstractname}{}

\begin{abstract}
\noindent \textbf{Abstract.} We propose a simple and original approach for solving  linear-quadratic mean-field stochastic control problems. We study both finite-horizon and infinite-horizon pro\-blems, and allow notably  some coefficients to be stochastic. Extension to the common noise case is also addressed. 
Our  method  is based on a suitable version  of the martingale formulation for verification theorems in control theory. The optimal control involves the solution to a system of Riccati ordinary differential equations and to a linear mean-field backward stochastic differential equation; existence and uniqueness conditions are provided for such a system. Finally, we illustrate our results through  an application to the production of an exhaustible resource. 
\end{abstract}
 
\vspace{3mm}


\vspace{3mm}


\section{Introduction}
\label{Sec:intro}

In recent years, optimal control of McKean-Vlasov stochastic differential equations, i.e., equations involving the law of the state process, has gained more and more attention, due to the increasing importance of problems with mean-field interactions and problems with cost functionals depending on the law of the state process and/or the law of the control (e.g., mean-variance portfolio selection problems or risk measures in finance). The goal of this paper is to design an elementary original approach for solving  linear-quadratic McKean-Vlasov control problems, which provides a unified framework for a wide class of problems and allows to treat problems which, to our knowledge, have not been studied before (e.g., common noise  and stochastic coefficients in the infinite-horizon case).

Linear-quadratic McKean-Vlasov (LQMKV) control problems are usually tackled by calculus of variations methods via stochastic maximum principle and decoupling techniques. Instead, we here consider a different approach, based on an extended version of the standard martingale formulation for verification theorems in control theory. Our approach  is valid for both finite-horizon and infinite-horizon problems and is closely connected to the dynamic programming principle (DPP), which holds in the LQ framework by taking into account both the state and its mean, hence restoring the time consistency of the problem. Notice that \cite{pha16} also used a DPP approach, but in the Wasserstein space of probability measures, and considered a priori closed-loop controls. Our approach is simpler in the sense that it does not rely on the notion of derivative in the Wasserstein space, and considers the larger class of  open-loop controls. We are able to obtain analytical solutions via the resolution of a system of two Riccati equations and the solution to a  linear mean-field backward stochastic differential equation. 

We first consider linear-quadratic McKean-Vlasov (LQMKV) control problems in finite horizon, where we allow some coefficients to be stochastic. We prove, by means of a weak martingale optimality principle, that there exists, under mild assumption on the coefficients, a unique optimal control, expressed in terms of the solution to a suitable system of Riccati equations and SDEs. We then provide some alternative sets of assumptions for the coefficient. We also show how the results adapt to the case where several independent Brownian motions are present. We also consider problem with common noise: Here, a similar formula holds, now considering conditional expectations. We then study the infinite-horizon case, characterizing the optimal control and the value function. Finally, we propose a detailed application, dealing with an infinite-horizon model of production of an exhaustible resource with a large number of producers and random price process.

We remark that in the infinite-horizon case some additional assumptions on the coefficients are required. On the one hand, having a well-defined value function requires a lower bound on the discounting parameter. On the other hand, we here deal with an infinite-horizon SDE, and the existence of a solution is a non-trivial problem. Finally, the admissibility of the optimal control requires a further condition of the discounting coefficient.

The literature on McKean-Vlasov control problems is now quite important, and we refer to the recent books by Bensoussan, Frehse and Yam \cite{benetal13} and Carmona and Delarue \cite{cardel18}, and the references therein. In this McKean-Vlasov framework, linear-quadratic (LQ) models provide an important class of solvable applications, and have been studied in many papers, including \cite{benetal11,Yong2013,HuangLiYong,pha16,gra16}, however  mostly for constant or deterministic coefficients, with the exception of  \cite{Sun} on finite horizon, 
and  \cite{LiSunYong}, which   deals with stochastic coefficients but considering a priori closed-loop strategies in linear form w.r.t.~the state and its mean.

The contributions of this paper are the following. 
First, we provide a new elementary solving technique for linear-quadratic McKean-Vlasov (LQMKV) control problems both on finite and infinite horizon. Second, the approach we propose has the advantage of being adaptable to several problems and allows several generalizations, which have not yet been studied before, as here outlined.   
In particular, we are able to solve  the case with common noise with some random coefficients, in finite and infinite horizon.  
The only references to this class of problems are the paper \cite{gra16} on finite horizon where the coefficients are deterministic,  and the paper \cite{pha16}, where the controls are required to be adapted to the filtration of the common noise (we here consider the case where $\alpha$ is adapted more generally to the pair of Brownian motions, that is, the one in the SDE and the common noise). As in  \cite{Sun}, we allow some coefficients to be stochastic, but to the best of our knowledge, this is the first time that explicit formulas are provided for infinite-horizon McKean-Vlasov control pro\-blems with random coefficients in the payoff. 
The inclusion of randomness in  some coefficients  is an important point, as it leads to a richer class of models, which is useful for many  applications, see e.g.~the investment problem in distributed generation under a random centralised electricity price studied in  \cite{AidBaseiPham}. 

The paper is organized as follows. Section \ref{Sec:finitepb} introduces finite-horizon LQMKV problems. Section \ref{Sec:assumpt}  presents the precise assumptions on the coefficients of the problems and provides a detailed description of the solving technique. In Section \ref{Sec:procedure} we solve, step by step, the control problem. Some remarks on the assumptions, and extensions are collected in Section \ref{Sec:remarks}. In Section \ref{Sec:infinitepb} we adapt the results to the infinite-horizon case. An application is studied  in Section \ref{Sec:applications}, which combines common noise and random coefficients.  
Finally, Section \ref{Sec:concl} concludes.

\section{Formulation of the Finite-Horizon Problem}
\label{Sec:finitepb}

Given a finite horizon $T>0$ (in Section \ref{Sec:infinitepb} we will extend the results to the infinite-horizon case), we fix a filtered probability space $(\Omega,\Fc,\F,\P)$, where $\F$ $=$ $(\Fc_t)_{0 \leq t \leq T}$  satisfies the usual conditions and is the natural filtration of a standard real Brownian motion $W$ $=$ $(W_t)_{0 \leq t \leq T}$, 
augmented with an independent $\sigma$-algebra $\Gc$.   Let $\rho \geq 0$ be a discount factor, and define the set of admissible (open-loop) controls as
\beqs
\aaa &:=&  \left\{ \alpha: \Omega \times [0,T] \to \rr^m \text{ s.t.~$\alpha$ is $\F$-adapted  and}\int_0^T e^{-\rho t} \eee[|\alpha_t|^2] dt < \infty \right\}.
\enqs
Given a square-integrable $\Gc$-measurable random variable $X_0$,  
and a control  $\alpha \in \aaa$, we consider the controlled linear 
mean-field stochastic differential equation in $\rr^d$ defined by
\begin{equation}
\label{pb:SDE}
\begin{cases}
dX^\alpha_t = b_t\big(X^\alpha_t, \eee[X^\alpha_t], \alpha_t, \eee[\alpha_t]\big) dt + \sigma_t\big(X^\alpha_t, \eee[X^\alpha_t], \alpha_t, \eee[\alpha_t]\big) dW_t, \quad 0 \leq t \leq T, 
\\
X^\alpha_0=X_0,
\end{cases}
\end{equation}
where for each $t \in [0,T]$, $x,\bar x \in \rr^d$ and $a,\bar a \in \rr^m$ we have set
\begin{equation}
\label{pb:coeffSDE}
\begin{array}{ccc}
b_t\big(x, \bar x, a, \bar a \big) & := &  \beta_t + B_t x + \tilde B_t \bar x + C_t a + \tilde C_t \bar a, \\
\sigma_t\big(x, \bar x, a, \bar a \big) & := &  \gamma_t + D_t x + \tilde D_t \bar x + F_t a + \tilde F_t \bar a.
\end{array}
\end{equation}
Here, the coefficients $\beta,\gamma$ of the affine terms are vector-valued $\F$-progressively measurable processes, whereas the other coe\-fficients of the linear terms are deterministic matrix-valued processes, see Section \ref{Sec:assumpt} for precise assumptions.  The quadratic cost functional to be minimized over $\alpha$ $\in$ $\Ac$ is 
\begin{equation} \label{pb:payoff}
\begin{array}{ccl}
J(\alpha) & := &  \E \Big[ \Int_0^T e^{-\rho t} f_t\big(X^\alpha_t, \eee[X^\alpha_t], \alpha_t, \eee[\alpha_t]\big) dt + e^{-\rho T}
g\big(X^\alpha_T, \eee[X^\alpha_T] \big) \Big], \\ \\
\rightarrow \;\;\;  V_0 & := &  \Inf_{\alpha \in \aaa} J(\alpha),
\end{array}
\end{equation}
where, for each $t \in [0,T]$, $x,\bar x \in \rr^d$ and $a,\bar a \in \rr^m$ we have set
\begin{equation} \label{pb:coeffPayoff}
\begin{array}{ccl}
f_t\big(x, \bar x, a, \bar a \big) &:=& (x-\bar x)\trans Q_t (x-\bar x) + \bar x\trans(Q_t + \tilde Q_t) \bar x +   2a\trans I_t (x-\bar x) + 2 \bar a\trans (I_t+\tilde I_t) \bar x   \\
& & \;\;\; + \;  (a - \bar a)\trans N_t (a - \bar a)  + \bar a\trans(N_t + \tilde N_t) \bar a +  2M_t\trans x + 2H_t\trans a, \\ 
g\big(x, \bar x\big) &:=&  (x - \bar x)\trans P (x - \bar x) + \bar x\trans( P+ \tilde P) \bar x + 2L\trans x.  
\end{array}
\end{equation}
Here, the coefficients $M,H,L$ of the linear terms are vector-valued $\F$-progressively measurable processes, whereas the other coefficients are deterministic matrix-valued processes. 
We refer again  to Section \ref{Sec:assumpt} for the precise assumptions. 
The symbol $\trans$ denotes the transpose of any vector or matrix. 

\vspace{1mm}

\begin{remark}
	{\rm {\bf a.} We have centred in \reff{pb:coeffPayoff}  the quadratic terms in the payoff functions $f$ and $g$. One could equivalently formulate the quadratic terms in non-centred form as    
		\begin{equation*} \label{pb:coeffPayoff2}
		\begin{array}{ccl}
		\tilde f_t\big(x, \bar x, a, \bar a \big) &:=& x\trans Q_t x + \bar x\trans \tilde Q_t \bar x  + a\trans N_t a 
		+ \bar a\trans \tilde N_t \bar a + 2M_t\trans x  + 2H_t\trans a + 2a\trans I_t x + 2 \bar a\trans \tilde I_t \bar x, \\ 
		\tilde g\big(x, \bar x\big) &:=&  x\trans P x + \bar x\trans \tilde P \bar x + 2L\trans x,
		\end{array}
		\end{equation*}
		by noting that, since $Q$, $P$,  $N$, $I$  are assumed to be deterministic, we have
		\begin{equation*}
		\E[ f_t\big(X^\alpha_t, \eee[X^\alpha_t], \alpha_t, \eee[\alpha_t]\big) ] = \E[ \tilde f_t\big(X^\alpha_t, \eee[X^\alpha_t], \alpha_t, \eee[\alpha_t]\big) ], \quad
		\E[g(X^\alpha_T, \E[X^\alpha_T])] = \E[\tilde g(X^\alpha_T, \E[X^\alpha_T])].
		\end{equation*}

		\vspace{1mm}
		
		\noindent {\bf b.} Notice that the only coefficients allowed to be stochastic are $\beta,\gamma,M,H,L$. Moreover, we note that in \eqref{pb:payoff}-\reff{pb:coeffPayoff} we could also consider a term of the form $\tilde M_t^T\bar x$, and then reduce for free the resulting problem to the case $\tilde M_t=0$. Indeed, since we consider the expectation of the running cost, we could equivalently substitute such a term with $\eee[\tilde M_t]\trans x$ by noting that 
		$\E[ \tilde M_t\trans \E[X_t^\alpha]]$ $=$ $\E[ \E[\tilde M_t]\trans X_t^\alpha]$. 
		Similarly, we do not need to consider  terms $\tilde H_t\trans \bar a$ and $\bar a\trans \check I_t x$, $a\trans \check I_t \bar x$ (for a deterministic matrix $\check I_t$). 
		
		\vspace{1mm}
		
		\noindent {\bf c.} See Section \ref{Sec:remarks} for the case where several Brownian motions and a common noise are present.
	}
	\qed 
\end{remark}

\section{Assumptions and Verification Theorem}
\label{Sec:assumpt}

Throughout the paper, for each $q \in \nn$ we denote by $\mathbb{S}^q$ the set of $q$-dimensional symmetric matrices. Moreover, for each normed space $(\M, |\cdot|)$ we set
\begin{align*}
L^{\infty}([0,T],\M) &:= \bigg\{ \phi : [0,T]\to \M \text{ s.t.~$\phi$ is measurable and $\textstyle \sup_{t \in [0,T]} |\phi_t| < \infty$} \bigg\}, 
\\
L^{2}([0,T],\M) &:= \bigg\{ \phi : [0,T]\to \M \text{ s.t.~$\phi$ is measurable and $\int_0^T e^{-\rho t} |\phi_t|^2 dt < \infty$} \bigg\},
\\
L^2_{\Fc_T}(\M) &:= \Big\{ \phi : \Omega \to \M \text{ s.t.~$\phi$ is $\Fc_T$-measurable and $\eee[|\phi|^2] < \infty$} \Big\}, 
\\
\Sc^2_{\F}(\Omega \times [0,T],\M) &:= 
\bigg\{ \phi : \Omega \times [0,T] \to \M \text{ s.t.~$\phi$ is $\F$-prog.~meas.~and $\E[\sup_{0\leq t\leq T} |\phi_t|^2] < \infty$} \bigg\},
\\
L^2_{\F}(\Omega \times [0,T],\M) &:= \bigg\{ \phi : \Omega \times [0,T]\to \M \text{ s.t.~$\phi$ is $\F$-progr.~meas.~and $\int_0^T e^{-\rho t} \eee[|\phi_t|^2] dt < \infty$} \bigg\}.
\end{align*}

We ask the following conditions on the coefficients of the problem to hold in the finite-horizon case.

\begin{itemize}
	\item[\textbf{(H1)}] The coefficients in \eqref{pb:coeffSDE} satisfy:
	\begin{itemize}
	 \item[(i)] 	$\beta,\gamma \in L^2_{\F}(\Omega \times [0,T],\rr^d)$,
	 \item[(ii)]  $B,\tilde B, D, \tilde D \in L^\infty([0,T],\rr^{d\times d})$,   $C,\tilde C, F, \tilde F \in L^\infty([0,T],\rr^{d\times m})$.
	 \end{itemize}
 	\item[\textbf{(H2)}] The coefficients in \eqref{pb:coeffPayoff} satisfy:
	\begin{itemize}
	\item[(i)] $Q, \tilde Q \in L^\infty([0,T],\mathbb{S}^d)$, $P, \tilde P \in \mathbb{S}^d$,  $N, \tilde N \in L^{\infty}([0,T],\mathbb{S}^m)$, $I, \tilde I \in L^\infty([0,T],\rr^{m\times d})$,
	\item[(ii)] $M \in L^2_{\F}(\Omega \times [0,T],\rr^d)$, $H \in L^2_{\F}(\Omega \times [0,T],\rr^m)$, $L \in L^2_{\mathcal{F}_T}(\rr^d)$,
	\item[(iii)] there exists $\delta > 0$ such that, for each $t \in [0,T]$,
	\begin{align} \nonumber 
	& N_t \; \geq \;  \delta \, \I_m,
	&& P \; \geq \;  0,
	&& Q_t \!-\! I_t\trans N^{-1}_t I_t \; \geq \;  0,
	\end{align}
	\item[(iv)] there exists $\delta > 0$ such that, for each $t \in [0,T]$,
	\begin{align} \nonumber
	& N_t \!+\! \tilde N_t \; \geq \;  \delta \, \I_m,
	&& P \!+\! \tilde P \; \geq \;  0,
	&& (Q_t \!+\! \tilde Q_t) \!-\! (I_t \!+\! \tilde I_t)\trans (N_t \!+\! \tilde N_t)^{-1} (I_t \!+\! \tilde I_t) \; \geq \;  0.
	\end{align}
	\end{itemize}
\end{itemize}

\vspace{1mm}
\begin{remark}
{\rm The uniform positive definite assumption on $N$ and $N+ \tilde N$ is a standard and natural coercive condition when dealing with linear-quadratic control problems. 
We discuss in Section 4 (see Remark \ref{remH2})  alternative assumptions when $N$ and $\tilde N$ may be degenerate. 
}
\qed
\end{remark}

\vspace{2mm}

By {\bf (H1)} and classical results, e.g.~\cite[Prop. 2.6]{Yong2013}, there exists a unique strong solution $X^\alpha = (X^\alpha_t)_{0 \leq t \leq T}$ to the mean-field SDE \eqref{pb:SDE}, which satisfies the standard estimate 
\begin{equation}
\label{pb:estimX}
\eee\big[ \sup_{t \in [0,T]}  |X^\alpha_t|^2\big] \leq C_\alpha \big (1+\eee[|X_0|^2] \big) < \infty,
\end{equation}
where $C_\alpha$ is a constant which depends on $\alpha$ $\in$ $\Ac$ only via $\int_0^T e^{-\rho t} \eee[|\alpha_t|^2] dt$. Also, by {\bf (H2)} and \eqref{pb:estimX}, the LQMKV control problem  \eqref{pb:payoff} is well-defined, in the sense that 
\begin{equation*}
\text{$J(\alpha) \in \rr$, for each $\alpha \in \aaa$.}
\end{equation*}

To solve the LQMKV control problem, we are going to use a suitable verification theorem. Namely, we consider an extended version of the martingale optimality principle usually cited in stochastic control theory: see Remark \ref{rem:verif} for a discussion.

\begin{Lemma}[\textbf{Finite-horizon verification theorem}]
\label{prop:verif}
Let $\{\Wc_t^\alpha, t \in [0,T], \alpha \in \aaa\}$ be a family of $\F$-adapted processes in the form $\Wc_t^\alpha$ $=$ 
$w_t(X^\alpha_t, \eee[X^\alpha_t])$ for some $\F$-adapted random field $\{w_t(x,\bar x), t\in [0,T], x,\bar x \in \R^d\}$ satisfying
\beq \label{growthw}
w_t(x,\bar x) & \leq & C( \chi_t + |x|^2 + |\bar x|^2), \;\;\; t \in [0,T], \; x,\bar x \in \R^d, 
\enq
for some positive constant $C$, and nonnegative process $\chi$ with $\sup_{t\in [0,T]} \E[|\chi_t|]$ $<$ $\infty$,  and such that
\begin{itemize}
\item[(i)] $w_T(x,\bar x)$ $=$  $g(x,\bar x)$, $x,\bar x$ $\in$ $\R^d$; 
\item[(ii)] the map $t$ $\in$ $[0,T]$ $\longmapsto$ $\E[\Sc_t^\alpha]$, with 
$\Sc_t^\alpha$ $:=$ $e^{-\rho t} \Wc_t^\alpha + \int_0^t e^{-\rho s}  f_s\big(X^\alpha_s, \eee[X^\alpha_s], \alpha_s, \eee[\alpha_s]\big) ds$, is nondecreasing for all 
$\alpha \in \aaa$;
\item[(iii)]  the map $t$ $\in$ $[0,T]$ $\longmapsto$ $\E[\Sc_t^{\alpha^*}]$ is constant for some $\alpha^* \in \aaa$.
	\end{itemize}
Then, $\alpha^*$ is an optimal control and $\eee[w_0(X_0,\eee[X_0])]$ is the value of the LQMKV control problem  \eqref{pb:payoff}: 
\begin{equation*}
V_0 \; = \;  \eee[ w_0(X_0,\eee[X_0])] \; = \;  J(\alpha^*).
\end{equation*}
Moreover, any other optimal control satisfies the condition (iii). 
\end{Lemma}
\noindent {\bf Proof.} From the growth condition \reff{growthw} and estimation \reff{pb:estimX}, we see that the function 
\beqs
t \in [0,T] & \longmapsto & \E [ \Sc_t^\alpha]  
\enqs
is well-defined for any $\alpha$ $\in$ $\Ac$. By (i) and (ii), we have for all $\alpha \in \aaa$ 
\beqs
\E[w_0(X_0,\E[X_0])] & = & \E[\Sc_0^\alpha] \\ 
& \leq & \E[\Sc_T^\alpha] \; = \;  \E\big[ e^{-\rho T} g(X_T^\alpha,\E[X_T^\alpha]) +  \int_0^T e^{-\rho t}  f_t\big(X^\alpha_t, \eee[X^\alpha_t], \alpha_t, \eee[\alpha_t]\big) dt \big] \\
& & \hspace{1.15cm}  = \; J(\alpha), 
\enqs
which shows that $\E[w_0(X_0,\E[X_0])]$ $\leq$ $V_0$ $=$ $\inf_{\alpha\in\Ac} J(\alpha)$, since $\alpha$ is arbitrary. Moreover, condition (iii)  with $\alpha^*$ shows that 
$\E[w_0(X_0,\E[X_0])]$ $=$ $J(\alpha^*)$, which proves the optimality of $\alpha^*$ with $J(\alpha^*)$ $=$ $\E[w_0(X_0,\E[X_0])]$.  Finally, suppose that $\tilde\alpha$ $\in$ $\Ac$ is another optimal control. Then
\beqs
\E[\Sc_0^{\tilde\alpha}] \; = \;  \E[w_0(X_0,\E[X_0])]  \; = \; J(\tilde\alpha) &=& \E[\Sc_T^{\tilde\alpha}].
\enqs
Since the map $t$ $\in$ $[0,T]$ $\longmapsto$ $\E[\Sc_t^{\tilde\alpha}]$ is nondecreasing, this shows that this map is actually constant, and concludes the proof. 
\qed
 
\vspace{3mm}

The general procedure to apply such a verification theorem consists of  the following three steps.
\begin{itemize}
\item[-] \emph{Step 1}. We guess a suitable parametric expression for the candidate random field $w_t(x,\bar x)$, 
and set for each $\alpha \in \aaa$ and $t \in [0,T]$,
	\begin{equation}
	\label{defS}
	\Sc^\alpha_t := e^{-\rho t} w_t(X^\alpha_t, \eee[X^\alpha_t]) + \int_0^t e^{-\rho s}  f_s\big(X^\alpha_s, \eee[X^\alpha_s], \alpha_s, \eee[\alpha_s]\big) ds.
	\end{equation}
	\item[-] \emph{Step 2}.  We apply It\^o's formula to $\Sc_t^\alpha$, for $\alpha$ $\in$ $\Ac$, and take the expectation to get
\beqs
d\E[\Sc_t^\alpha] &=& e^{-\rho t} \E[ \Dc_t^\alpha ] dt, 
\enqs
for some $\F$-adapted processes $\Dc^\alpha$ with 
\beqs
 \E[ \Dc_t^\alpha]  &=&  \E \Big[ - \rho w_t(X_t^\alpha,\E[X_t^\alpha]) + \frac{d}{dt} \E\big[ w_t(X_t^\alpha,\E[X_t^\alpha])\big] + 
f_t(X_t^\alpha,\E[X_t^\alpha],\alpha_t,\E[\alpha_t]) \Big]. \nonumber 
\enqs
We then determine  the coefficients  in the random field $w_t(x,\bar x)$  s.t.~condition (i) in Lemma \ref{prop:verif} (i.e., $w_T(.)$ $=$ $g(.)$) is satisfied, and so as to have
\beqs \label{conDc}
 \; \E[\Dc_t^\alpha] \; \geq \; 0,  \; t \geq 0, \forall \alpha\in \Ac, &\mbox{and}&  \; \E[\Dc_t^{\alpha^*}] \; = \; 0, \; t \geq 0, 
\mbox{ for some } \;  \alpha^* \in \Ac,
\enqs
which ensures that the mean optimality principle conditions (ii) and (iii) are satisfied, and then $\alpha^*$ will be the optimal control. This leads to a system of backward ordinary and stochastic differential equations.
 \item[-] \emph{Step 3}. We study the existence of solutions to the system obtained in Step 2, which will also ensure the square integrability condition of $\alpha^*$ in $\Ac$, hence its optimality.  
\end{itemize}

\vspace{1mm}
\begin{remark} \label{rem:verif}
{\rm
The standard martingale optimality principle used in the verification theorem for stochastic control problems, see e.g.~\cite{elk81}, consists in finding 
a family of processes $\{\Wc_t^\alpha, 0\leq t\leq T, \alpha\in\Ac\}$ s.t. 
\begin{itemize}
\item[(ii')] the process $\Sc_t^\alpha$ $=$ 
$e^{-\rho t} \Wc_t^\alpha + \int_0^t e^{-\rho s}  f_s\big(X^\alpha_s, \eee[X^\alpha_s], \alpha_s, \eee[\alpha_s]\big) ds$, $0\leq t\leq T$, is a submartingale for each $\alpha$ $\in$ 
$\Ac$, 
\item[(iii')]  $\Sc^{\alpha^*}$ is a martingale for some $\alpha^*$ $\in$ $\Ac$, 
\end{itemize}
which obviously implies the weaker conditions (ii) and (iii)  in Proposition \ref{prop:verif}. Practically, the martingale optimality conditions (ii')-(iii')  would reduce via the It\^o decomposition of $\Sc^\alpha$ to the condition that $\mathcal{D}_t^\alpha$ $\geq$ $0$, for each $\alpha \in \aaa$, and $\mathcal{D}_t^{\alpha^*}$ $=$ $0$, $0\leq t\leq T$,  for a suitable control  $\alpha^*$. In the classical framework of stochastic control problem without mean-field dependence, one looks for $\Wc_t^\alpha$ $=$ $w_t(X_t^\alpha)$ for some random field $w_t(x)$ depen\-ding only on the state value, and the martingale optimality principle leads to the classical Hamilton-Jacobi-Bellman (HJB) equation (when all the coefficients are non-random) or to a stochastic HJB, see \cite{pen92}, in the general random coefficients case. In our context of McKean-Vlasov control problems, one looks for $\Wc_t^\alpha$ $=$ $w_t(X_t^\alpha,\E[X_t^\alpha])$ depending also on the mean of the state value, and the pathwise condition 
on $\Dc_t^\alpha$ would not allow us  to determine a suitable random field $w_t(x,\bar x)$. Instead, we exploit the weaker condition (ii) formulated as a mean condition on $\E[\Dc_t^\alpha]$, and we shall see in the next section how it leads indeed to a suitable characterization of  $w_t(x,\bar x)$. The methodology of the weak martingale approach in Lemma 3.2 works concretely whenever one can find a family of value functions for the McKean-Vlasov control problem that depends upon the law of the state process only via its mean (or conditional mean in the case of common noise). This imposes a Markov property on the pair of controlled process $(X_t,\mathbb{E}[X_t])$, and hence a linear structure of the dynamics for $X$ w.r.t. its mean and the control. The running payoff and terminal cost function $f,g$ should then also  depend on the state and its mean (or conditional mean in the case of common noise), but not necessarily in the quadratic form. The quadratic form has the advantage of suggesting a suitable quadratic form for the candidate value function, while in general it is not explicit. 	Actually, the candidate $w_t(x,x')$  for the value function should satisfy a Bellman PDE in finite dimension, namely the dimension of $(X_t,\mathbb{E}[X_t])$, which is a particular finite dimensional case of the Master equation. This argument  of making the McKean-Vlasov control problem finite-dimensional  is exploited more generally in \cite{BalataEtAl} where the dependence on the law is through the first $p$-moments of the state process.
}
\qed
\end{remark}

\section{Solution to LQMKV}

\label{Sec:procedure}

In this section, we apply the verification theorem in Lemma \ref{prop:verif} to characterize an optimal control for the problem  \eqref{pb:payoff}. We will follow the procedure outlined at the end of Section \ref{Sec:assumpt}. In the sequel, and for convenience of notations, we set
\begin{equation} \label{nothat}
\begin{aligned}{}
&\hat B_t \; := \; B_t + \tilde B_t, \quad\,\,\, &&\hat C_t \; := \; C_t + \tilde C_t,  \quad\,\,\, &&\hat D_t \; := \; D_t + \tilde D_t, \quad\,\,\, &&\hat F_t \; := \; F_t + \tilde F_t, \\
&\hat I_t \; := \; I_t + \tilde I_t, \quad\,\,\, &&\hat N_t \; := \; N_t + \tilde N_t,  \quad\,\,\, &&\hat Q_t \; := \; Q_t + \tilde Q_t, \quad\,\,\, &&\hat P \; := \; P  + \tilde P. 
\end{aligned}
\end{equation}

\vspace{1mm}
\begin{remark}	\label{rem:notations}
{\rm To simplify the notations, throughout the paper we will often denote expectations by an upper bar and omit the dependence on controls. Hence, for example, we will write $X_t$ for $X^\alpha_t$, $\bar X_t$ for $\eee[X^\alpha_t]$, and $\bar\alpha_t$ for $\E[\alpha_t]$. 
}
\qed
\end{remark}

\noindent {\bf Step 1.}  Given the quadratic structure of the cost functional $f_t$ in \eqref{pb:coeffPayoff}, we infer a candidate for the random field  $\{w_t(x,\bar x), t\in [0,T], x,\bar x \in \R^d\}$ in the form:
\beq \label{wquadra}
w_t(x,\bar x) &=& (x - \bar x)\trans K_t (x - \bar x) + \bar x\trans \Lambda_t \bar x + 2 Y_t\trans x + R_t,
\enq
where $K_t,\Lambda_t,Y_t,R_t$ are  suitable processes to be determined later. 
The centering of the quadratic term is a convenient choice, which provides simpler calculations. Actually, since the quadratic coefficients in the payoff \reff{pb:coeffPayoff} are deterministic symmetric matrices, we  look for deterministic symmetric matrices $K,\Lambda$ as well. Moreover, since in the statement of Lemma \ref{prop:verif} we always consider the expectation of $\Wc^\alpha_t$ $=$ $w_t(X_t^\alpha,\E[X_t^\alpha])$, we can assume, w.l.o.g., that $R$ is deterministic.  Given the randomness of the linear coefficients in \reff{pb:coeffPayoff},  
the process $Y$ is considered in general as an $\F$-adapted process. 
Finally, the terminal condition $w_T(x,\bar x)$ $=$ $g(x,\bar x)$ ((i) in Lemma \ref{prop:verif}) determines the terminal conditions satisfied by $K_T,\Lambda_T,Y_T,R_T$. We then search for processes  $(K,\Lambda,Y,R)$ valued in $\S^d\times \S^d\times \R^d\times \R$ in backward form: 
\begin{equation}\label{KLYfirst}
\left\{
\begin{array}{cclcl}
dK_t &=& \dot K_t dt, & 0 \leq t \leq T, & K_T \; = \; P, \\
d\Lambda_t & = & \dot \Lambda_t dt,  &  0 \leq t \leq T, & \Lambda_T \; = \;  \hat P, \\
dY_t & = & \dot Y_t dt + Z^Y_t dW_t, & 0 \leq t \leq T, & Y_T \; = \; L, \\
dR_t &=& \dot R_t dt, & 0 \leq t \leq T, & R_T \; = \;  0,
\end{array}
\right.
\end{equation}
for some deterministic processes $(\dot K,\dot\Lambda,\dot R)$ valued in $\S^d\times\S^d\times\R$, and $\F$-adapted processes $\dot Y, Z^Y$ valued in $\R^d$. 

\vspace{3mm}

\noindent {\bf Step 2.} For $\alpha \in \aaa$ and $t \in [0,T]$, let $\Sc^\alpha_t$ as in \reff{defS}. We have 
\beqs
d\E[\Sc_t^\alpha] &=& e^{-\rho t} \E[ \Dc_t^\alpha ] dt, 
\enqs
for some $\F$-adapted processes $\Dc^\alpha$ with 
\beqs
\E[ \Dc_t^\alpha] &=&  \E \big[ - \rho w_t(X_t^\alpha,\E[X_t^\alpha]) + \frac{d}{dt} \E\big[ w_t(X_t^\alpha,\E[X_t^\alpha])\big] + 
f_t(X_t^\alpha,\E[X_t^\alpha],\alpha_t,\E[\alpha_t]) \big]. \nonumber 
\enqs
We apply the It\^o's formula to $w_t(X_t^\alpha,\E[X_t^\alpha])$, recalling the quadratic form \reff{wquadra} of $w_t$, the equations in \reff{KLYfirst}, and the dynamics (see equation \reff{pb:SDE})
\beqs
d \bar X_t^\alpha &=& [\bar\beta_t + \hat B_t \bar X_t^\alpha + \hat C_t \bar\alpha_t] dt,  \\
d(X_t^\alpha- \bar X_t^\alpha) & = & \big[ \beta_t - \bar\beta_t + B_t (X_t^\alpha- \bar X_t^\alpha) + C_t (\alpha_t - \bar\alpha_t) \big] dt \\
& &  + \; \big[  \gamma_t + D_t(X_t^\alpha-\bar X_t^\alpha) + \hat D_t \bar X_t^\alpha + F_t(\alpha_t-\bar\alpha_t) + \hat F_t \bar\alpha_t \big] dW_t,
\enqs
where we use the upper bar notation for expectation, see Remark \ref{rem:notations}. Recalling the quadratic form \reff{pb:coeffPayoff} of the running cost $f_t$, we obtain, after careful but straightforward computations, that
\beq 
\E[ \Dc_t^\alpha] 
&=&  \E \Big[ (X_t - \bar X_t)\trans (\dot K_t + \Phi_t ) (X_t - \bar X_t) 
+ \bar X_t\trans\big(\dot\Lambda_t + \Psi_t \big) \bar X_t  \nonumber \\
& & \hspace{2cm}   + \; 2\big( \dot  Y_t + \Delta_t  \big)\trans X_t  + \dot R_t - \rho R_t  + \bar \Gamma_t   +    \chi_t(\alpha)  \Big],    \label{expressD} 
\enq
(we omit the dependence in $\alpha$ of $X$ $=$ $X^\alpha$, $\bar X$ $=$ $\bar X^\alpha$), where
\begin{equation} \label{PhiK}
\left\{
\begin{array}{rcl}
\Phi_t &:=& - \rho K_t+ K_tB_t +  B_t\trans K_t + D_t\trans K_t D_t + Q_t \; = \; \Phi_t(K_t), \\
\Psi_t &:=& - \rho \Lambda_t +  \Lambda_t\hat B_t + \hat B_t\trans\Lambda_t + \hat D_t\trans K_t \hat D_t + \hat Q_t \; = \; \Psi_t(K_t,\Lambda_t), \\
\Delta_t&:=&  - \rho Y_t + B_t\trans Y_t + \tilde B_t\trans\bar Y_t + D_t\trans Z_t^Y + \tilde D_t\trans\overline{Z_t^Y} + K_t (\beta_t - \bar\beta_t) + \Lambda_t \bar\beta_t \\
& & \;\;\;\;\;  + \;  D_t\trans K_t (\gamma_t - \bar \gamma_t) + \hat D_t\trans K_t \bar\gamma_t + M_t \; = \; \Delta_t(K_t,\Lambda_t,Y_t,\bar Y_t,Z_t^Y,\overline{Z_t^Y}), \\
\Gamma_t &:=&   \gamma_t\trans K_t \gamma_t +  2 \beta_t\trans Y_t + 2 \gamma_t\trans Z_t^Y \; = \; \Gamma_t(K_t,Y_t,Z_t^Y),
\end{array}
\right.
\end{equation}
for $t$ $\in$ $[0,T]$, and
\beq 
\chi_t(\alpha) &:=&  (\alpha_t - \bar \alpha_t)\trans S_t (\alpha_t - \bar\alpha_t) + \bar\alpha_t\trans \hat S_t  \bar\alpha_t \nonumber  \\
& & \;\;\; + \;  2 \big(  U_t (X_t-\bar X_t)  + V_t \bar X_t  + O_t + \xi_t - \bar\xi_t  \big)\trans \alpha_t.  \label{defchi}
\enq
Here,  the deterministic coefficients $S_t, \hat S_t, U_t, V_t, O_t$ are defined, for $t$ $\in$ $[0,T]$, by
\begin{equation} \label{defSUVO}
\left\{
\begin{array}{rcl}
S_t&:=&  N_t + F_t\trans K_t F_t \;=\; S_t(K_t), \\
\hat S_t &:=& \hat N_t + \hat F_t\trans K_t \hat F_t \;=\; \hat S_t(K_t), \\
U_t &:=& I_t + F_t\trans K_t D_t + C_t\trans K_t \;=\; U_t(K_t), \\
V_t &:=& \hat I_t + \hat F_t\trans K_t \hat D_t + \hat C_t\trans \Lambda_t \;=\; V_t(K_t,\Lambda_t), \\
O_t &:=& \bar H_t + \hat F_t\trans K_t \bar \gamma_t + \hat C_t\trans \bar Y_t + \hat F_t\trans \overline{Z_t^Y} \;=\; O_t(K_t,\bar Y_t, \overline{Z_t^Y}) ,
\end{array}
\right.
\end{equation}
and the stochastic coefficient $\xi_t$ 
of mean $\bar\xi_t$ is 
defined, for $t$ $\in$ $[0,T]$, by
\begin{equation} \label{defxi}
\left\{
\begin{array}{ccl}
\xi_t & := & H_t + F_t\trans K_t \gamma_t  + C_t\trans Y_t + F_t\trans Z_t^Y \;=\; \xi_t(K_t,Y_t,Z_t^Y), \\
\bar\xi_t & := & \bar H_t + F_t\trans K_t \bar \gamma_t  + C_t\trans \bar Y_t + F_t\trans \overline{Z_t^Y} \;=\; \bar\xi_t(K_t,\bar Y_t,\overline{Z_t^Y}).  
\end{array}
\right.
\end{equation}
Notice that we have suitably rearranged the terms in \reff{expressD} in order to keep only linear terms in $X$ and $\alpha$, by using the elementary observation that $\E[ \phi_t\trans \bar X_t]$ $=$ $\E[\bar \phi_t \trans X_t]$, and $\E[\psi_t\trans\bar \alpha_t]$ $=$ $\E[\bar \psi_t\trans\alpha_t]$ for any vector-valued random variable $\phi_t$, $\psi_t$ of mean $\bar\phi_t$, $\bar\psi_t$.   
 
Next, the key-point is to complete the square w.r.t.~the control $\alpha$  in the process $\chi_t(\alpha)$ defined in  \reff{defchi}. Assuming for the moment that the symmetric matrices $S_t$ and $\hat S_t$ are positive definite in $\S^m$ (this will follow typically from the non-negativity of the matrix $K$, as checked in Step 3, and conditions (iii)-(iv) in {\bf (H2)}), it is clear that  one can find  a   deterministic $\R^{m\times m}$-valued $\Theta$   (which may be not unique) s.t.  for all $t$ $\in$ $[0,T]$, 
\beqs
\Theta_t S_t \Theta_t\trans &=& \hat S_t,  
\enqs
for all $t$ $\in$ $[0,T]$, and which is also deterministic like $S_t,\hat S_t$. We can then rewrite the expectation of $\chi_t(\alpha)$ as  
\beqs
\E[\chi_t(\alpha)] &=& \E \Big[ \big( \alpha_t  - \bar\alpha_t + \Theta_t\trans\bar\alpha_t - \eta_t)\trans S_t \big( \alpha_t  - \bar\alpha_t + \Theta_t\trans\bar\alpha_t - \eta_t \big) - \zeta_t \Big],
\enqs
where 
\beqs
\eta_t &:=& a_t^0(X_t,\bar X_t) + \Theta_t\trans a_t^1(\bar X_t),
\enqs
with $a_t^0(X_t,\bar X_t)$ a centred random variable, and $a_t^1(\bar X_t)$ a deterministic function
\beqs
a_t^0(x,\bar x) \; := \;   - S_t^{-1} U_t (x-\bar x) -  S_t^{-1} (\xi_t-\bar \xi_t),  & & a_t^1(\bar x) \; := \;   -   \hat S_t^{-1} (V_t \bar x + O_t),  
\enqs
and
\beqs
\zeta_t &:=&  (X_t-\bar X_t) \trans\big( U_t\trans S_t^{-1} U_t \big)(X_t-\bar X_t) + \bar X_t\trans \big(V_t\trans \hat S_t^{-1}V_t\big) \bar X_t \\
& & \;\;\; + \;   2\big(U_t\trans S_t^{-1}(\xi_t-\bar \xi_t)  + V_t\trans \hat S_t^{-1} O_t \big)\trans X_t  \\
& & \;\;\; + \; (\xi_t - \bar\xi_t)\trans S_t^{-1} (\xi_t - \bar\xi_t)    +   O_t\trans \hat S_t^{-1} O_t.
\enqs 
We can then  rewrite the expectation in \reff{expressD} as 
\beqs 
 & & \E[\Dc_t^\alpha] \\
&=&  \E \Big[ (X_t - \bar X_t)\trans \big(\dot K_t + \Phi_t^0 \big) (X_t - \bar X_t)  + \;  \bar X_t\trans\big(\dot\Lambda_t + \Psi_t^0 \big) \bar X_t  \nonumber \\
& & \;   + \; 2\big( \dot  Y_t + \Delta_t^0 \big)\trans X_t   + \; \dot R_t - \rho R_t  + \overline{\Gamma_t^0}    \\
& & \; + \;   \big( \alpha_t - a_t^0(X_t,\bar X_t)  - \bar\alpha_t + \Theta_t\trans(\bar\alpha_t - a_t^1(\bar X_t))\big)\trans S_t \big( \alpha_t - a_t^0(X_t,\bar X_t)  - \bar\alpha_t + \Theta_t\trans(\bar\alpha_t - a_t^1(\bar X_t)) \big)  \Big],
\enqs
where we set 
\begin{equation} \label{Phi0}
\left\{
\begin{array}{rcl}
\Phi_t^0 &:=& \Phi_t -   U_t\trans S_t^{-1} U_t \; = \; \Phi_t^0(K_t), \\
\Psi_t^0 &:=& \Psi_t  -   V_t\trans \hat S_t^{-1}V_t \; = \; \Psi_t^0(K_t,\Lambda_t), \\
\Delta_t^0 &:=& \Delta_t -  U_t\trans S_t^{-1}(\xi_t-\bar \xi_t)  - V_t\trans \hat S_t^{-1} O_t \;= \; \Delta_t^0(K_t,\Lambda_t,Y_t,\bar Y_t,Z_t^Y,\overline{Z_t^Y}), \\
\Gamma_t^0 &:= & \Gamma_t - (\xi_t - \bar\xi_t)\trans S_t^{-1} (\xi_t - \bar\xi_t)    -   O_t\trans \hat S_t^{-1} O_t \; = \;  
\Gamma_t^0(K_t,Y_t,\bar Y_t,Z_t^Y,\overline{Z_t^Y}),
\end{array}
\right.
\end{equation}
and stress the dependence on $(K,\Lambda,Y,Z^Y)$ in view of \reff{PhiK}, \reff{defSUVO}, \reff{defxi}.   
Therefore, whenever 
\beqs
\dot K_t + \Phi_t^0  \; = \;  0, \qquad \dot\Lambda_t + \Psi_t^0  \;=\;  0, & &  
\dot  Y_t + \Delta_t^0 \; = \;  0, \qquad  \dot R_t - \rho R_t  + \overline{\Gamma_t^0} \; = \;  0
\enqs
holds for all $t$ $\in$ $[0,T]$, we have
\beq
 & & \E[\Dc_t^\alpha] \label{Dcalphafin} \\
&=& \E \Big[  \big( \alpha_t - a_t^0(X_t,\bar X_t)  - \bar\alpha_t + \Theta_t\trans(\bar\alpha_t - a_t^1(\bar X_t))\big)\trans S_t \big( \alpha_t - a_t^0(X_t,\bar X_t)  - \bar\alpha_t + \Theta_t\trans(\bar\alpha_t - a_t^1(\bar X_t)) \big)   \Big], \nonumber
\enq
which is non-negative for all $0\leq t\leq T$, $\alpha$ $\in$ $\Ac$, i.e.,  the process $\Sc^\alpha$ satisfies the condition (ii) of the verification theorem in Lemma \ref{prop:verif}. We are then led to consider the  following system of backward (ordinary and stochastic) differential equations (ODEs and BSDE):
\begin{equation} \label{sysK}
\left\{
\begin{array}{ccl}
dK_t &=&  -  \Phi_t^0(K_t) dt, \;\;\; 0 \leq t \leq T, \; K_T \; = \; P, \\ 
d\Lambda_t &=& -  \Psi_t^0(K_t,\Lambda_t)   dt,  \;\;\; 0 \leq t \leq T, \; \Lambda_T \; = \; P + \tilde P, \\
dY_t &=& - \Delta_t^0(K_t,\Lambda_t,Y_t,\E[Y_t],Z_t^Y,\E[Z_t^Y]) dt + Z_t^Y dW_t, \;\; 0 \leq t \leq T, \; Y_T \; = \; L,  \\
dR_t &=& \big[ \rho R_t - \E[\Gamma_t^0(K_t,Y_t, \E[Y_t],Z_t^Y,\E[Z_t^Y]) \big] dt, \;\;\; 0 \leq t \leq T, \; R_T \; = \; 0. 
\end{array}
\right.
\end{equation}

\begin{Definition}
A solution to the system \reff{sysK} is a quintuple of processes $(K,\Lambda,Y,Z^Y,R)$ $\in$ 
$L^\infty([0,T],\S^d)\times L^\infty([0,T],\S^d)\times \Sc^2_{\F}(\Omega\times [0,T],\R^d)\times L^2_{\F}(\Omega\times [0,T],\R^d)\times L^\infty([0,T],\R)$ 
s.t.~the $\S^m$-valued  processes $S(K)$, $\hat S(K)$ $\in$ $L^\infty([0,T],\S^m)$, are positive definite a.s., and the following relation 
\begin{equation*}
\left\{
\begin{array}{ccl}
K_t &=& P + \int_t^T \Phi_s^0(K_s) ds, \\
\Lambda_t &=& P + \tilde P + \int_t^T \Psi_s^0(K_s,\Lambda_s) ds, \\
Y_t &=& L + \int_t^T \Delta_s^0(K_s,\Lambda_s,Y_s,\E[Y_s],Z_s^Y,\E[Z_s^Y]) ds +  \int_t^T Z_s^Y dW_s, \\
R_t &=& \int_t^T \big( - \rho R_s + \E\big[\Gamma_s^0(K_s,Y_s, \E[Y_s],Z_s^Y,\E[Z_s^Y]) \big] \big) ds,
\end{array}
\right.
\end{equation*}
holds for all $t$ $\in$ $[0,T]$. 
\end{Definition}

We shall discuss in the next paragraph (Step 3) the existence of a solution to the system of ODEs-BSDE \reff{sysK}. For the moment, we provide the connection between this system and the solution to the LQMKV control problem. 

\begin{Proposition} \label{profini}
Suppose that $(K,\Lambda,Y,Z^Y,R)$ is a solution to the system of ODEs-BSDE \reff{sysK}. Then, the control process
\beqs
\alpha_t^* &=& a_t^0(X_t^*,\E[X_t^*]) + a_t^1(\E[X_t^*]) \\
&=& - S_t^{-1}(K_t) U_t(K_t) (X_t-\E[X_t^*]) -  S_t^{-1}(K_t) \big(\xi_t(K_t,Y_t,Z_t^Y) - \bar \xi_t(K_t,\E[Y_t],\E[Z_t^Y]\big) \\
& & \;\;\;  -  \;  \hat S_t^{-1}(K_t) \big(V_t(K_t,\Lambda_t) \E[X_t^*] +  O_t(K_t,\E[Y_t], \E[Z_t^Y]) \big),
\enqs
where $X^*$ $=$ $X^{\alpha^*}$ is the state process with the feedback control $a_t^0(X_t^*,\E[X_t^*]) + a_t^1(\E[X_t^*])$,  is the optimal control for the LQMKV problem \reff{pb:payoff}, i.e., 
$V_0$ $=$ $J(\alpha^*)$, and we have
\beqs
V_0 &=& \E\big[(X_0 - \E[X_0])\trans K_0 (X_0 - \E[X_0]) \big]  + \E[X_0]\trans \Lambda_0 \E[X_0] + 2 \E[Y_0\trans X_0]  + R_0. 
\enqs 
\end{Proposition}
\noindent {\bf Proof.} Consider a solution $(K,\Lambda,Y,Z^Y,R)$ to the system \reff{sysK}, and let $w_t$ as of the quadratic form \reff{wquadra}. First, notice that $w$ satisfies the growth condition \reff{growthw} as $K,\Lambda,R$ are bounded and $Y$ satisfies a square-integrability condition  in $L^2_{\F}(\Omega\times [0,T],\R^d)$. The terminal condition $w_T(.)$ $=$ $g$ is also satisfied from the terminal condition of the system \reff{sysK}. Next, for this choice of $(K,\Lambda,Y,Z^Y,R)$, the expectation $\E[\Dc_t^\alpha]$ in \reff{Dcalphafin}  is non\-negative for all $t$ $\in$ $[0,T]$, $\alpha$ $\in$ $\Ac$, which means that the process $\Sc^\alpha$  satisfies the condition (ii) of the verification theorem in Lemma \ref{prop:verif}. Moreover, we see that $\E[\Dc_t^{\alpha}]$ $=$ $0$, $0\leq t\leq T$, for some $\alpha$ $=$ $\alpha^*$ if and only if (recall that $S_t$ is positive definite a.s.)
\beqs 
\alpha_t^* - a_t^0(X_t^*,\E[X_t^*])  - \E[\alpha_t^*] + \Theta_t\trans(\E[\alpha_t^*] - a_t^1(\E[X_t^*])) &=& 0, \;\;\; 0 \leq t \leq T.  
\enqs
Taking expectation in the above relation, and recalling that $\E[a_t^0(X_t^*,\E[X_t^*])]$ $=$ $0$, $\Theta_t$ is invertible, we get $\E[\alpha_t^*]$ $=$ $a_t^1(\E[X_t^*])$, and thus 
\beq \label{expressalpha}
\alpha_t^* &=& a_t^0(X_t^*,\E[X_t^*]) + a_t^1(\E[X_t^*]), \;\;\; 0 \leq t \leq T.  
\enq
Notice that $X^*$ $=$ $X^{\alpha^*}$ is solution to a linear McKean-Vlasov dynamics, and satisfies the square-integrability condition $\E[\sup_{0\leq t\leq T}|X_t^*|^2]$ $<$ $\infty$,  which implies in its turn that $\alpha^*$ satisfies the square-integrability condition  $L^2_{\F}(\Omega\times [0,T],\R^m)$, since $S^{-1}$, $\hat S^{-1}$, $U$, $V$ are bounded, 
and $O$, $\xi$ are square-integrable respectively in  $L^2([0,T],\R^m)$ and $L_{\F}^2(\Omega\times [0,T],\R^m)$. Therefore, $\alpha^*$ $\in$ $\Ac$, and we conclude by the verification 
theorem in Lemma \ref{prop:verif} that it is the unique optimal control.
\qed

\vspace{2mm}

\paragraph{Step 3.} Let us now verify under assumptions {\bf (H1)}-{\bf (H2)} the existence and uniqueness of a  solution to the decoupled system in \eqref{sysK}.

\begin{itemize}
	\item[(i)] We first consider the equation for $K$, which is actually a matrix Riccati equation written as:
\begin{equation} \label{eqK}
\left\{
\begin{array}{rcl}
\frac{d}{dt} K_t + Q_t -\rho K_t + K_tB_t+ B_t\trans K_t + D_t\trans K_t D_t & & \\
- (I_t + F_t\trans K_t D_t + C_t\trans K_t)\trans(N_t + F_t\trans K_t F_t)^{-1}(I_t + F_t\trans K_t D_t + C_t\trans K_t) &=& 0, \; t \in [0,T], \\
K_T &=& P.
\end{array}
\right.
\end{equation}
Multi-dimensional Riccati equations are known to be related to control theory. Namely,  \eqref{eqK} is associated  to the standard linear-quadratic stochastic control problem:
\beqs
v_t(x) &:=& \inf_{\alpha\in \Ac} \E \Big[ \int_t^T e^{-\rho t} \Big( (\tilde X_s^{t,x,\alpha})\trans Q_s \tilde X_s^{t,x,\alpha} + 2 \alpha_s\trans I_s \tilde X_s^{t,x,\alpha} 
+ \alpha_s\trans N_s\alpha_s \Big) ds \\
& & \hspace{2cm} + \; e^{-\rho T}  (\tilde X_T^{t,x,\alpha})\trans P  \tilde X_T^{t,x,\alpha}  \Big] ,
\enqs
where $\tilde X^{t,x,\alpha}$ is the controlled linear dynamics solution to
\beqs
d\tilde X_s &=& ( B_s \tilde X_s + C_s \alpha_s) ds + (D_s \tilde X_s + F_s \alpha_s) dW_s, \;\;\; t \leq s \leq T, \; \tilde X_t = x. 
\enqs
 By a standard result in control theory (see \cite[Ch.~6, Thm.~6.1, 7,1, 7.2]{YongZhou}, with a straightforward adaptation of the arguments to include the discount factor), 
 under {\bf (H1)}, {\bf (H2)}(i)-(ii),  there exists a unique solution $K$ $\in$ $L^\infty([0,T],\S^d)$ with  $K_t$ $\geq$ $0$  to \eqref{eqK}, provided that
	\begin{equation}
	\label{ipoK}
	P \; \geq \;  0, \qquad Q_t - I_t\trans N_t^{-1} I_t \; \geq \; 0, \qquad N_t \; \geq \;  \delta \, \I_m, \;\;\; 0 \leq t \leq T, 
	\end{equation}
	for some  $\delta>0$, which is true by  {\bf (H2)}(iii), and in this case, we have $v_t(x)$ $=$ $x\trans K_t x$.  Notice also that $S(K)$ $=$ $N+ F\trans KF$ is positive definite. 
	
	\item[(ii)] Given $K$, we now consider the equation for $\Lambda$. Again, this is a matrix  Riccati equation that we rewrite as
	\begin{equation} \label{eqL}
\left\{
\begin{array}{rcl}
\frac{d}{dt} \Lambda_t + \hat Q^K_t -\rho \Lambda_t + \Lambda_t\hat B_t+ \hat B_t \trans  \Lambda_t & & \\
- \big(\hat I^K_t + \hat C_t\trans\Lambda_t \big)\trans(\hat N^K_t)^{-1} \big(\hat I^K_t + \hat C_t\trans \Lambda_t \big)  &=& 0, \; t \in [0,T], \\
\Lambda_T &=& \hat P ,
\end{array}
\right.
\end{equation}
where we have set, for $t \in [0,T]$,
\beqs
\hat Q_t^K &:=& \hat Q_t   + \hat D_t \trans K_t \hat D_t, \\
\hat I^K_t &:=&  \hat I_t + \hat F_t \trans K_t \hat D_t, \\
\hat N^K_t &:=&   \hat N_t + \hat F_t  \trans K_t \hat F_t.
\enqs
As for the equation for $K$, there exists a unique solution $\Lambda$ $\in$ $L^\infty([0,T],\S^d)$ with  $\Lambda_t$ $\geq$ $0$  to \eqref{eqL}, provided that
	\begin{equation}
	\label{ipoLambda}
	\hat P \; \geq \;  0, \;\;\; \hat Q^K_t - (\hat I^K_t)\trans (\hat N^K_t)^{-1} (\hat I^K_t) \; \geq \;  0, \;\;\;  \hat N^K_t \; \geq \;  \delta \, \I_m, \;\;\; 0 \leq t \leq T, 
	\end{equation}
	for some  $\delta>0$.  
Let us check that {\bf (H2)}(iv) implies \eqref{ipoLambda}. We already have $\hat P$ $\geq$ $0$.  Moreover, as $K\geq 0$ we have: 
$\hat N^K_t$ $\geq$ $\hat N_t$ $\geq$ $\delta \I_m$.  By simple algebraic manipulations and as $\hat N_t > 0$, we have (omitting the time dependence)
\begin{align*}
	\hat Q^{K} - (\hat I^{K})\trans(\hat N^{K})^{-1} \hat I^{K} & = \hat Q - \hat I\trans \hat N^{-1} \hat I + (\hat D - \hat F \hat N^{-1} \hat I)\trans K (\hat D - \hat F \hat N^{-1} \hat I) 
	\\
	& \qquad - \Big(\hat F\trans K (\hat D - \hat F \hat N^{-1} \hat I) \Big)\trans(\hat N + \hat F\trans K \hat F)^{-1} \Big(\hat F \trans K (\hat D - \hat F \hat N^{-1} \hat I) \Big)
	\\
	& \geq \;  \hat Q - \hat I\trans \hat N^{-1} \hat I + (\hat D - \hat F \hat N^{-1} \hat I)\trans K (\hat D - \hat F \hat N^{-1} \hat I) 
	\\
	& \qquad - \Big(\hat F\trans K (\hat D - \hat F \hat N^{-1} \hat I) \Big)\trans (\hat F\trans K \hat F)^{-1} \Big(\hat F\trans K (\hat D - \hat F \hat N^{-1} \hat I) \Big)
	\\
	& = \; \hat Q - \hat I\trans \hat N^{-1} \hat I  \; \geq \;  0, \;\;\; \mbox{ by {\bf (H2)}(iv)}.
\end{align*} 		
	\item[(iii)] Given $(K,\Lambda)$, we consider the equation for $(Y,Z^Y)$. This is a mean-field linear BSDE written as
 	\begin{equation}
	\label{eqY}
	\begin{cases}
	dY_t = \Big( \vartheta_t + G_t\trans(Y_t - \eee[Y_t]) + \hat G_t\trans \eee[Y_t] + J_t\trans(Z_t^Y - \eee[Z_t^Y]) 
	+ \hat{J}_t\trans\eee[Z_t^Y] \Big) dt + Z_t^Y dW_t, 
	\\
	Y_T = L,
	\end{cases}
	\end{equation}
	where the deterministic coefficients $G$, $\hat G$, $J$, $\hat J$ $\in$ $L^\infty([0,T],\R^{d\times d})$, and the stochastic process  $\vartheta$  $\in$ $L^2_{\F}(\Omega\times [0,T],\R^d)$ 
	are defined by
	\begin{align*}
	& G_t \; := \;  \rho \,\, \I_d  - B_t + C_tS^{-1}_t U_t, 
	\\
	& \hat G_t \; := \;  \rho \,\, \I_d - \hat B_t + \hat C_t\hat S^{-1}_t V_t,
	\\
	& J_t \; := \;  - D_t + F_tS^{-1}_t U_t,
	\\
	& \hat{J}_t \; := \;  - \hat D_t + \hat F_t  \hat S^{-1}_t V_t, 
	\\
	&\vartheta_t \; := \;  - M_t - K_t (\beta_t - \eee[\beta_t] ) - \Lambda_t \eee[\beta_t] - D_t\trans K_t (\gamma_t - \eee[\gamma_t]) - \hat D_t\trans K_t \eee[\gamma_t]
	\\
	& \qquad  \;\;\; + \;  U_t\trans S^{-1}_t\big( H_t - \eee[H_t] + F_t\trans K_t (\gamma_t - \eee[\gamma_t])\big) \;+\;  
	 V_t\trans \hat S^{-1}_t \big(\eee[H_t] + \hat F_t\trans K_t \eee[\gamma_t] \big),
	\end{align*}
and the expressions for $S,\hat S,U,V$ are recalled in \eqref{defSUVO}. By standard results, see \cite[Thm.~2.1]{LiSunXiong}, there exists a unique solution $(Y,Z^Y)$ $\in$ 
$\Sc_{\F}^2(\Omega\times [0,T],\R^d)\times L_{\F}^2(\Omega\times [0,T],\R^d)$ to \eqref{eqY}.  
	\item[(iv)]  Given $(K,\Lambda,Y,Z^Y)$, the equation for $R$ is a linear ODE, whose unique solution is explicitly given  by 
	\beq \label{eqR}
	R_t &=& \int_t^T e^{-\rho (s-t)} h_s ds.
	\enq
	Here, the deterministic function $h$ is defined, for $t \in [0,T]$, by
	\beqs
	h_t &:=&  \eee\big[ - \gamma_t\trans K_t \gamma_t - \beta_t\trans Y_t - 2\gamma_t\trans Z_t^Y  + \xi_t\trans S^{-1}_t \xi_t \big] 
	- \eee[\xi_t]\trans S^{-1}_t \eee[\xi_t] + O_t\trans \hat S^{-1}_t O_t,
	\enqs
	and the expressions of $O$ and $\xi$ are recalled in \reff{defSUVO} and \reff{defxi}. 
\end{itemize}

\vspace{2mm}
To sum up the arguments of this section, we have proved the following result.

\begin{Theorem}
	\label{thm:optimal}
	Under assumptions  {\bf (H1)}-{\bf (H2)}, the optimal control for the LQMKV  pro\-blem  \eqref{pb:payoff} is given  by
	\beq
		\label{optimal}
		\alpha^*_t &=&  - S_t^{-1} U_t(X^*_t - \eee[X^*_t]) - \hat S_t^{-1} (V_t \eee[X^*_t]+ O_t) - S^{-1}_t (\xi_t-\eee[\xi_t]),
	\enq
	where $X^*=X^{\alpha^*}$ and the deterministic coefficients $S,\hat S$  $\in$ $L^\infty([0,T],\S^m)$, $U,V$ in $L^\infty([0,T],\R^{m\times d})$, $O$ $\in$  
	$L^\infty([0,T],\R^{m})$ and the stochastic coefficient $\xi$ $\in$ $L^2_{\F}(\Omega\times [0,T],\R^m)$ are defined by
	\begin{equation}
		\label{coeffTHM}
		\left\{
		\begin{array}{ccl}
			S_t &:=& N_t + F_t\trans K_t F_t,
			\\
			\hat S_t &:=& N_t + \tilde N_t + (F_t + \tilde F_t)\trans K_t (F_t + \tilde F_t),
			\\
			U_t & :=& I_t + F_t\trans K_t D_t + C_t\trans K_t,
			\\
			V_t & :=& I_t+\tilde I_t + (F_t + \tilde F_t)\trans K_t (D_t + \tilde D_t) + (C_t+\tilde C_t)\trans \Lambda_t,
			\\
			O_t & :=& \eee[H_t] + (F_t + \tilde F_t)\trans K_t \eee[\gamma_t] + (C_t + \tilde C_t)\trans  \eee[Y_t] + (F_t + \tilde F_t)\trans \eee[Z_t^Y],
			\\
			\xi_t & :=&  H_t + F_t\trans K_t \gamma_t + C_t\trans Y_t + F_t\trans Z_t^Y,
		\end{array}
		\right. 
	\end{equation}
	with  $(K,\Lambda,Y,Z^Y,R)$ $\in$ 
	$L^\infty([0,T],\S^d)\times L^\infty([0,T],\S^d)\times \Sc^2_{\F}(\Omega\times [0,T],\R^d)\times L^2_{\F}(\Omega\times [0,T],\R^d)\times L^\infty([0,T],\R)$ the unique  solution to \reff{eqK}, \reff{eqL}, \reff{eqY}, \reff{eqR}. The corresponding value of the problem is 
	\beqs
		V_0 &=& J(\alpha^*) \; = \;  \eee\big[ (X_0- \eee[X_0])\trans K_0 (X_0- \eee[X_0]) \big] +  \eee[X_0]\trans \Lambda_0  \eee[X_0] + 2 \eee\big[Y_0\trans X_0] + R_0.
	\enqs
\end{Theorem}

\section{Remarks and Extensions}
\label{Sec:remarks}

We collect here some remarks and extensions for the problem presented in the previous sections.

\vspace{1mm}
\begin{remark} \label{remH2}
{\rm  Assumptions {\bf (H2)}(iii)-(iv) are used only for ensuring the existence of a non-negative solution $(K,\Lambda)$ to equations \reff{eqK}, \reff{eqL}. In some specific cases, they can be substituted by alternative conditions. 
	
For example, in the one-dimensional case $n$ $=$ $m$ $=$ $1$ (real-valued control and state variable), with $N$ $=$ $0$ (no quadratic term on the control in the running cost) and $I$ $=$ $0$, the equation for $K$ writes
\begin{equation*}
\frac{d}{dt}K_t + Q_t + \big( -\rho + 2B_t - C^2_t/F^2_t -2C_tD_t/F_t \big) K_t = 0, \,\,\, t \in [0,T], \qquad K_T = P.
\end{equation*}
This is a first-order linear ODE, which clearly admits a unique solution, provided that $F_t\neq 0$. Moreover, when $P$ $>$ $0$, then $K$ $>$ $0$ by classical maximum principle, so that we have $S_t$ $>$ $0$. Hence, an alternative condition to {\bf (H2)}(iii) is, for $t \in [0,T]$,

\vspace{4mm}

\noindent {\bf (H2)}(iii') \hspace{1.7cm}  $n$ $=$ $m$ $=$ $1$, $N_t$ $=$ $I_t$ $=$ $0$, $P$ $>$ $0$, $F_t$ $\neq$ $0$. 

\vspace{4mm}

Let us now discuss an alternative condition to the  uniform positive condition on $N+\tilde N$ in {\bf (H2)}(iv), in the case where $N+\tilde N$ is only assumed to be non-negative. 
When the constant matrix $P$ is positive definite, then $K$ is uniformly positive definite in $\S^d$, i.e., $K_t$ $\geq$ $\delta \I_m$, $0\leq t\leq T$, for some $\delta$ $>$ $0$, by strong maximum principle for the ODE \reff{eqK}. Then, when $F+\tilde F$ is  uniformly non-degenerate, i.e., $|F_t+\tilde F_t|$ $\geq$ $\delta$, $0\leq t\leq T$, for some $\delta$ $>$ $0$, we see that $\hat S_t$ $=$ $\hat N_t^K$ $\geq$ $(F_t+\tilde F_t)\trans K_t (F_t+\tilde F_t)$ $\geq$ $\delta' \I_d$ for some $\delta'$ $>$ $0$. Notice also that when $I+\tilde I$ $=$ $0$, then $\hat Q^K - (\hat I^K)\trans (\hat N^K)^{-1} (\hat I^K)$ $\geq$ $Q+\tilde Q$.  Consequently, assumption {\bf (H2)}(iv) can be alternatively replaced by 

\vspace{4mm}

\noindent {\bf (H2)}(iv') \hspace{1cm} $N_t+\tilde N_t$,  $P+\tilde P$, $Q_t+\tilde Q_t$ $\geq$ $0$,  $P$ $>$ $0$, $I_t+\tilde I_t$ $=$ $0$,  $|F_t+\tilde F_t|$ $\geq$ $\delta$,  

\vspace{4mm}

\noindent for $t \in [0,T]$ and some $\delta$ $>$ $0$, which ensures that condition \reff{ipoLambda} is satisfied, hence giving  the existence and uniqueness of a nonnegative solution $\Lambda$ to \reff{eqL}. 

We underline that {\bf (H2)}(iii')-(iv') are not the unique alternative to {\bf (H2)}(iii)-(iv). In some applications, none of such conditions is satisfied, typically as $Q$ $=$ $\tilde Q$ $=$ $0$, while $I$ or $\tilde I$ is non-zero. However, a solution $(K,\Lambda)$ (possibly non-positive) to \reff{eqK}-\reff{eqL} may still exist, with $S(K)$ and $\hat S(K)$ positive definite, and one can then still apply Proposition \ref{profini} to get the conclusion of Theorem \ref{thm:optimal}, i.e., the optimal control exists and is given by \reff{optimal}.
}
\qed
\end{remark} 


\vspace{1mm}

\begin{remark} \label{remWmulti}
{\rm
The result in Theorem \ref{thm:optimal} can be easily extended to the case where several Brownian motions are present in the controlled equation:
\beqs
dX^\alpha_t &=& b_t\big(X^\alpha_t, \eee[X^\alpha_t], \alpha_t, \eee[\alpha_t]\big) dt 
+ \Sum_{i=1}^n\sigma^i_t\big(X^\alpha_t, \eee[X^\alpha_t], \alpha_t, \eee[\alpha_t]\big) dW^i_t,
\enqs
where $W^1, \dots, W^n$ are standard independent real Brownian motions and, for each $t \in [0,T]$, $i\in\{1,\dots,n\}$, $x,\bar x \in \rr^d$ and $a,\bar a \in \rr^m$, we  set
	\begin{equation}
	\label{newsde}
	\begin{gathered}
	b_t\big(x, \bar x, a, \bar a \big) := \beta_t + B_t x + \tilde B_t \bar x + C_t a + \tilde C_t \bar a,
	\\
	\sigma^i_t\big(x, \bar x, a, \bar a \big) := \gamma^i_t + D^i_t x + \tilde D^i_t \bar x + F^i_t a + \tilde F^i_t \bar a.
	\end{gathered}
	\end{equation}
We ask the coefficients in \eqref{newsde} to satisfy a suitable adaptation of {\bf (H1)}: namely, we substitute $D,\tilde D, F, \tilde F$ with $D^i,\tilde D^i, F^i, \tilde F^i$, for $i \in \{1,\dots,n\}$. The cost functional and {\bf (H2)}  are unchanged.
	
The statement of Theorem \ref{thm:optimal} is easily adapted to this extended framework. To simplify the notations we use Einstein convention: for example, we write 
$(D^i_t)\trans K D^i_t$ instead of $\sum_{i=1}^n (D^i_t)\trans K D^i_t$. The optimal control $\alpha^*$ is given by \eqref{optimal}, where the coefficients are now defined by  
\begin{equation*}
		\left\{
		\begin{array}{ccl}
			S_t &:=& N_t + (F_t^i)\trans K_t F_t^i,
			\\
			\hat S_t &:=& N_t + \tilde N_t + (F_t^i + \tilde F_t^i)\trans K_t (F_t^i + \tilde F_t^i),
			\\
			U_t & :=& I_t + (F_t^i)\trans K_t D_t^i + C_t\trans K_t,
			\\
			V_t & :=& I_t+\tilde I_t + (F_t^i + \tilde F_t^i)\trans K_t (D_t^i + \tilde D_t^i) + (C_t+\tilde C_t)\trans \Lambda_t,
			\\
			O_t & :=& \eee[H_t] + (F_t^i + \tilde F_t^i)\trans K_t \eee[\gamma_t] + (C_t + \tilde C_t)\trans  \eee[Y_t] + (F_t^i + \tilde F_t^i)\trans \eee[Z_t^Y],
			\\
			\xi_t & :=&  H_t + (F_t^i)\trans K_t \gamma_t + C_t\trans Y_t + (F_t^i)\trans Z_t^Y,
		\end{array}
		\right. 
	\end{equation*}
and $(K,\Lambda,Y,Z^Y,R)$ $\in$ $L^\infty([0,T],\S^d)\times L^\infty([0,T],\S^d)\times \Sc^2_{\F}(\Omega\times [0,T],\R^d)\times L^2_{\F}(\Omega\times [0,T],\R^d)\times L^\infty([0,T],\R)$ is the unique  solution to \reff{sysK}, with 
\begin{equation} \nonumber
\left\{
\begin{array}{rcl}
\Phi_t(K_t) &=& - \rho K_t + K_tB_t +  B_t\trans K_t + (D_t^i)\trans K_t D_t^i + Q_t ,\\
\Psi_t(K_t,\Lambda_t) &=& - \rho \Lambda_t +  \Lambda_t(B_t+\tilde B_t) + (B_t+\tilde B_t)\trans\Lambda_t,  \\
& & \;\;\; + \;  (D_t^i+\tilde D_t^i)\trans K_t (D_t^i + \tilde D_t^i) + Q_t + \tilde Q_t, \\
\Delta_t(K_t,\Lambda_t,Y_t,\E[Y_t],Z_t^Y,\E[Z_t^Y]) &=&  - \rho Y_t + B_t\trans Y_t  + \tilde B_t\trans\E[Y_t] +  
k (\beta_t - \bar\beta_t) + \Lambda_t \bar\beta_t   + M_t  \\
& & \;   + \;  (D_t^i)\trans K_t \gamma_t + (\tilde D_t^i)\trans K_t \bar\gamma_t + (D_t^i)\trans Z_t^Y + (\tilde D_t^i)\trans\E[Z_t^Y], \\
\Gamma_t(K_t,Y_t,Z_t^Y) &=&   \gamma_t\trans K_t \gamma_t +  2 \beta_t\trans Y_t + 2 \gamma_t\trans Z_t^Y, \;\;\; t \in [0,T]. \hfill\text{\qed}
\end{array}
\right.
\end{equation}
}
\end{remark}

\vspace{1mm}

\begin{remark}
{\rm The optimal control provided by Theorem \ref{thm:optimal} generalizes known results and standard formulas in control theory. 
\begin{itemize}
\item[-] For example, in the case where
	\begin{equation*}
	I_t=\tilde I_t=\beta_t=\gamma_t=M_t=H_t=L_t=0,
	\end{equation*}
then $Y$ $=$ $Z^Y$ $=$ $0$, $R$ $=$ $0$ (correspondingly, we have $O$ $=$ $\xi$ $=$ $0$).  We thus retrieve the formula in \cite[Thm. 4.1]{Yong2013} for the 
optimal control (recalling the notations in \reff{nothat}):
\beqs
\alpha_t^* &=& - \big(N_t + F_t\trans K_t F_t\big)^{-1}(F_t\trans K_t D_t + C_t\trans K_t)(X^*_t - \eee[X^*_t]) \\
& & \; - \; \big(\hat N_t + \hat F_t  \trans K_t \hat F_t \big)^{-1} \big( \hat F_t\trans K_t \hat D_t + \hat C_t\trans \Lambda_t \big) \eee[X^*_t].
\enqs
\item[-] Consider now the case where all the mean-field coefficients are zero, that is 
	\begin{equation*}
	\label{no-mkv}
	\tilde B_t = \tilde C_t = \tilde D_t = \tilde F_t = \tilde Q_t = \tilde N_t = \tilde P_t \equiv 0.
	\end{equation*}
	Assume, in addition, that $\beta_t=\gamma_t=H_t=M_t=0$. In this case, $K_t$ $=$ $\Lambda_t$ satisfy the same Riccati equation, $Y_t=\hat Y_t=R_t=0$, and we have
	\beqs
	S_t &=& \hat S_t \; = \;  N_t + F_t\trans K_t F_t,	\\
	U_t &=&  V_t \; = \;  I_t + F_t\trans K_t D_t + C_t\trans K_t,	\\
	O &=& \xi \; = \;  0,
	\enqs
	which leads to the well-known formula for classical linear-quadratic control problems (see, e.g.~\cite{YongZhou}):
	\begin{equation*}
	\pushQED{\qed} 
	\alpha^*_t \; = \; - (N_t + F_t\trans K_t F_t)^{-1} ( I_t + F_t\trans K_t D_t + C_t\trans K_t)X^*_t, \;\;\; 0 \leq t\leq T.  \qedhere
	\popQED
	\end{equation*}
\end{itemize} 
}
\end{remark}

\vspace{1mm}

\begin{remark}
{\rm The mean  of the optimal state $X^*$ $=$ $X^{\alpha^*}$ can be computed as the solution of a linear ODE. Indeed,by plugging \eqref{optimal} into \eqref{pb:SDE} and taking expectation, we get
\beqs
\frac{d}{dt} \eee[X^*_t] &=&  \big( B_t+\tilde B_t - (C_t + \tilde C_t) \hat S_t^{-1}V_t \big) \eee[X^*_t] + \big( \eee[\beta_t] - (C_t + \tilde C_t) \hat S_t^{-1} O_t \big),
\enqs
which can be solved  explicitly in the one-dimensional case $d$ $=$ $1$, and expressed as an exponential of matrices in the multidimensional case. 
}
\qed
\end{remark}

\vspace{1mm}

\begin{remark} \label{remcommonnoise}
{\rm
\emph{(The case of common noise)}.  We now extend the results in Theorem \ref{thm:optimal} to the case where a common noise is present. Let $W$ and $W^0$ be two independent real Brownian motions defined on the same filtered probability space $(\Omega, \mathcal{F}_T,\mathbbm{F},\mathbbm{P})$. Let $\mathbb{F} = \{\mathcal{F}_t\}_{t \in [0,T]}$ be the filtration generated by the pair $(W,W^0)$ and let $\mathbb{F}^0 = \{\mathcal{F}^0_t\}_{t \in [0,T]}$ be the filtration generated by $W^0$.

For any $X_0$ and $\alpha \in \aaa$ as in Section \ref{Sec:intro}, the controlled process $X^\alpha_t$ is defined by
\begin{equation} \label{dynXnoise}
\begin{cases}
dX^\alpha_t = b_t\big(X^\alpha_t, \eee[X^\alpha_t|W^0_t], \alpha_t, \eee[\alpha_t|W^0_t]\big) dt + \sigma_t\big(X^\alpha_t, \eee[X^\alpha_t|W^0_t], \alpha_t, \eee[\alpha_t|W^0_t]\big) dW_t 
\\
\hspace{5cm} + \sigma^0_t\big(X^\alpha_t, \eee[X^\alpha_t|W^0_t], \alpha_t, \eee[\alpha_t|W^0_t]\big) dW^0_t, \quad 0 \leq t \leq T, 
\\
X^\alpha_0=X_0,
\end{cases}
\end{equation}
where for each $t \in [0,T]$, $x,\bar x \in \rr^d$ and $a,\bar a \in \rr^m$ we have set
\begin{equation*}
\begin{array}{ccc}
b_t\big(x, \bar x, a, \bar a \big) & := &  \beta_t + B_t x + \tilde B_t \bar x + C_t a + \tilde C_t \bar a, \\
\sigma_t\big(x, \bar x, a, \bar a \big) & := &  \gamma_t + D_t x + \tilde D_t \bar x + F_t a + \tilde F_t \bar a,
\\
\sigma^0_t\big(x, \bar x, a, \bar a \big) & := &  \,\,\,\, \gamma^0_t + D^0_t x + \tilde D^0_t \bar x + F^0_t a + \tilde F^0_t \bar a.
\end{array}
\end{equation*}
Here, $B$,$\tilde B$,$C$,$\tilde C$,$D$,$\tilde D$,$F$,$\tilde F$,$D^0$,$\tilde D^0$,$F^0$,$\tilde F^0$ are essentially bounded $\mathbbm{F}^0$-adapted processes, whereas $\beta, \gamma, \gamma^0$ are square-integrable $\mathbbm{F}$-adapted processes. We underline that $\beta, \gamma,\gamma^0$ can depend on $W$ as well. The problem is 
\begin{equation*} 
\begin{array}{ccl}
J(\alpha) & := &  \E \Big[ \Int_0^T e^{-\rho t} f_t\big(X^\alpha_t, \eee[X^\alpha_t|W^0_t], \alpha_t, \eee[\alpha_t|W^0_t]\big) dt + e^{-\rho T}
g\big(X^\alpha_T, \eee[X^\alpha_T|W^0_T] \big) \Big], \\
\rightarrow \;\;\;  V_0 & := &  \Inf_{\alpha \in \aaa} J(\alpha),
\end{array}
\end{equation*}
with $f_t, g$ as in \eqref{pb:coeffPayoff}. The coefficients in $f_t, g$ here satisfy the following assumptions: $Q$,$\tilde Q$,$I$,$\tilde I$,$N$,$\tilde N$ are essentially bounded $\mathbbm{F}^0$-adapted processes, $P ,\tilde P$ are essentially bounded $\mathcal{F}^0_T$-measurable random variables, $M,H$ are square-integrable $\mathbbm{F}$-adapted processes, $L$ is a square-integrable $\mathcal{F}_T$-measurable random variables. We also ask conditions (iii) and (iv) in \textbf{(H2)} to hold. We remark that $M,H,L$ can also depend on $W$.

As in Section \ref{Sec:procedure}, we guess a quadratic expression for the candidate random field. Namely, we consider $w_t(X^\alpha_t,\eee[X^\alpha_t|W^0_t])$, with $\{w_t(x,\bar x), t\in [0,T], x,\bar x \in \R^d\}$ as in \eqref{wquadra}, that we here recall:
\beq
w_t(x,\bar x) &=& (x - \bar x)\trans K_t (x - \bar x) + \bar x\trans \Lambda_t \bar x + 2 Y_t\trans x + R_t,
\enq
for suitable coefficients $K,\Lambda,Y,R$. Since the quadratic coefficients in $f_t,g$ are $\mathbbm{F}^0$-adapted, we guess that the coefficients $K,\Lambda$ are $\mathbbm{F}^0$-adapted as well (notice the difference with respect to Section \ref{Sec:procedure}, where $K,\Lambda$ were deterministic). The affine coefficients in $b_t,\sigma_t,\sigma^0_t$ and the linear coefficients in $f_t,g$ are $\mathbbm{F}$-adapted, so that $Y$ needs to depend on both $W$ and $W^0$. Finally, as in Section \ref{Sec:procedure} we can choose $R$ deterministic. We then look for processes $(K,\Lambda,Y,R)$ valued in $\S^d\times \S^d\times \R^d\times \R$ and in the form: 
\begin{equation*} 
\left\{
\begin{array}{cclcl}
dK_t &=& \dot K_t dt + Z^K_t dW^0_t, & 0 \leq t \leq T, & K_T \; = \; P, \\
d\Lambda_t & = & \dot \Lambda_t dt + Z^\Lambda_t dW^0_t,  &  0 \leq t \leq T, & \Lambda_T \; = \;  \hat P, \\
dY_t & = & \dot Y_t dt + Z^Y_t dW_t + Z^{Y,0}_t dW^0_t, & 0 \leq t \leq T, & Y_T \; = \; L, \\
dR_t &=& \dot R_t dt, & 0 \leq t \leq T, & R_T \; = \;  0,
\end{array}
\right.
\end{equation*}
for some $\F^0$-adapted processes $\dot K, \dot\Lambda, Z^K,Z^\Lambda$ valued in $\S^d$, some $\F$-adapted processes $\dot Y,Z^Y,Z^{Y,0}$ valued in $\R^d$ and a continuous function $\dot R$ valued in $\R$. 

We use the notations in \eqref{nothat} and extend them to the new coefficients $D^0,\tilde D^0,F^0,\tilde F^0$ (e.g., we denote $\hat D^0_t=D^0_t + \tilde D^0_t$). Moreover, for any random variable $\zeta$, we denote by $\bar \zeta$ the conditional expectation with respect ot $W^0_t$, i.e., $\bar \zeta = \eee[\zeta |W^0_t]$. For each $\alpha \in \aaa$ and $t \in [0,T]$, let $\Sc^\alpha_t$ be defined by
\begin{equation*}
\Sc^\alpha_t  \; := \;  e^{-\rho t} w_t(X^\alpha_t, \eee[X^\alpha_t|W^0_t]) + \int_0^t e^{-\rho s}  f_s\big(X^\alpha_s, \eee[X^\alpha_s|W^0_s], \alpha_s, \eee[\alpha_s|W^0_s]\big) ds,
\end{equation*}
and let $\Dc_t^\alpha$ be defined by $d\E[\Sc_t^\alpha] = e^{-\rho t} \E[ \Dc_t^\alpha ] dt$. By applying the It\^o formula to $\Sc^\alpha_t$, an expression for $\E[ \Dc_t^\alpha]$ is given by \eqref{expressD} and \eqref{defchi}, whose coefficients are now defined as follows. The coefficients in \eqref{expressD} are here given by
\begin{equation*} 
\left\{
\begin{array}{rcl}
\Phi_t &:=& - \rho K_t+ K_tB_t + B_t\trans K_t + Z^K_t D^0_t + (D^0_t)\trans Z^K_t 
\\
& & \;\;\;\;\; + D_t\trans K_t D_t + (D^0_t)\trans K_t D^0_t + Q_t \; = \; \Phi_t(K_t,Z^K_t), 
\\
\Psi_t &:=& - \rho \Lambda_t +  \Lambda_t\hat B_t + \hat B_t\trans\Lambda_t + Z^\Lambda_t\hat D^0_t + (\hat D^0_t)\trans Z^\Lambda_t 
\\
& & \;\;\;\;\; + \hat D_t\trans K_t \hat D_t + (\hat D^0_t)\trans \Lambda_t \hat D^0_t + \hat Q_t \; = \; \Psi_t(K_t,\Lambda_t,Z^\Lambda_t), 
\\
\Delta_t&:=&  - \rho Y_t + B_t\trans Y_t  + \tilde B_t\trans \overline{Y_t} + D_t\trans Z_t^Y + \tilde D_t\trans\overline{Z_t^Y} + (D^0_t)\trans Z^{Y,0}_t + (\tilde D^0_t)\trans\overline{Z_t^{Y,0}} 
\\
& & \;\;\;\;\; + K_t (\beta_t - \bar\beta_t) + \Lambda_t \bar\beta_t + Z^K_t (\gamma^0_t - \bar \gamma^0_t) + Z^\Lambda_t \bar \gamma^0_t + D_t\trans K_t (\gamma_t - \bar \gamma_t) + \hat D_t\trans K_t \bar\gamma_t 
\\
& & \;\;\;\;\;  + (D^0_t)' K_t (\gamma^0_t - \bar \gamma^0_t) + (\hat D^0_t)'\Lambda_t\bar \gamma^0_t + M_t \; = \; \Delta_t(K_t,\Lambda_t,Y_t,Z_t^Y, Z_t^{Y,0},\overline{Y_t},\overline{Z_t^Y},\overline{Z_t^{Y,0}}), 
\\
\Gamma_t &:=&   (\gamma_t - \bar \gamma_t)\trans K_t (\gamma_t - \bar \gamma_t) + \bar \gamma_t\trans K_t \bar \gamma_t + (\gamma^0_t - \bar \gamma^0_t)\trans K_t (\gamma^0_t - \bar \gamma^0_t) + (\bar \gamma^0_t)\trans \Lambda_t \bar \gamma^0_t 
\\
& & \;\;\;\;\; + 2 \beta_t\trans Y_t + 2 \gamma_t\trans Z_t^Y + 2 (\gamma^0_t)\trans Z_t^{Y,0}  \; = \; \Gamma_t(K_t,\Lambda_t, Y_t,Z_t^Y,Z_t^{Y,0}),
\end{array}
\right.
\end{equation*}
whereas the coefficients in \eqref{defchi} are given by
\begin{equation*} 
\left\{
\begin{array}{rcl}
S_t&:=&  N_t + F_t\trans K_t F_t + (F^0_t)\trans K_t F^0_t \;=\; S_t(K_t), \\
\hat S_t &:=& \hat N_t + \hat F_t\trans K_t \hat F_t +  (\hat F^0_t)\trans \Lambda_t \hat F^0_t \;=\; \hat S_t(K_t,\Lambda_t), \\
U_t &:=& I_t + F_t\trans K_t D_t + (F^0_t)\trans K_t D^0_t + C_t\trans K_t + (F^0_t)\trans Z^K_t \;=\; U_t(K_t,Z^K_t), \\
V_t &:=& \hat I_t + \hat F_t\trans K_t \hat D_t + (\hat F^0_t)\trans \Lambda_t \hat D^0_t + \hat C_t\trans \Lambda_t + (\hat F^0_t)\trans Z^\Lambda_t \;=\; V_t(K_t,\Lambda_t,Z^\Lambda_t), \\
O_t &:=& \bar H_t + \hat F_t\trans K_t \bar \gamma_t + (\hat F^0_t)\trans \Lambda_t \bar \gamma^0_t + \hat C_t\trans \overline{Y_t} + \hat F_t\trans \overline{Z_t^Y} + (\hat F^0_t)\trans \overline{Z^{Y,0}_t} \;=\; O_t(K_t,\overline{Y_t}, \overline{Z_t^Y},\overline{Z^{Y,0}_t}), \\
\xi_t & := & H_t + F_t\trans K_t \gamma_t + (F^0_t)\trans K_t \gamma^0_t + C_t\trans Y_t + F_t\trans Z_t^Y + (F^0_t)\trans Z_t^{Y,0} \;=\; \xi_t(K_t,Y_t,Z_t^Y,Z_t^{Y,0}).
\end{array}
\right.
\end{equation*}
Completing the square as in Section \ref{Sec:procedure} and setting to zero all the terms which do not depend on the control, we get that $(K,\Lambda,Y,R)$ satisfy the following problem
\begin{equation} \label{EDSRcommon}
\left\{
\begin{array}{ccl}
dK_t &:=&  -  \Phi_t^0 dt + Z^K_t dW^0_t, \qquad 0 \leq t \leq T, \qquad K_T \; = \; P, \\ 
d\Lambda_t &:=& -  \Psi_t^0 dt  + Z^\Lambda_t dW^0_t, \qquad 0 \leq t \leq T, \qquad \Lambda_T \; = \; P + \tilde P, \\
dY_t &:=& - \Delta_t^0 dt + Z_t^Y dW_t  + Z^{Y,0}_t dW^0_t, \qquad 0 \leq t \leq T, \quad Y_T \; = \; L,  \\
dR_t &:=& ( \rho R_t - \E[\Gamma_t^0]) dt, \qquad 0 \leq t \leq T, \qquad R_T \; = \; 0, 
\end{array}
\right.
\end{equation}
where the coefficients $\Phi^0,\Psi^0,\Delta^0,\Gamma^0$ are defined by
\begin{equation*} 
\left\{
\begin{array}{rcl}
\Phi_t^0 &:=& \Phi_t -   U_t\trans S_t^{-1} U_t \; = \; \Phi_t^0(K_t,Z^K_t), \\
\Psi_t^0 &:=& \Psi_t  -   V_t\trans \hat S_t^{-1}V_t \; = \; \Psi_t^0(K_t,\Lambda_t,Z^\Lambda_t), \\
\Delta_t^0 &:=& \Delta_t -  U_t\trans S_t^{-1}(\xi_t-\bar \xi_t)  - V_t\trans \hat S_t^{-1} O_t \;= \; \Delta_t^0(K_t, Z^K_t, \Lambda_t, Z^\Lambda_t, Y_t, Z_t^Y, Z_t^{Y,0}, \overline{Y_t}, \overline{Z_t^Y}, \overline{Z_t^{Y,0}}), \\
\Gamma_t^0 &:= & \Gamma_t - (\xi_t - \bar\xi_t)\trans S_t^{-1} (\xi_t - \bar\xi_t) - O_t\trans \hat S_t^{-1} O_t \; = \;  
\Gamma_t^0(Y_t, Z_t^Y, Z_t^{Y,0}, \overline{Y_t}, \overline{Z_t^Y}, \overline{Z_t^{Y,0}}).
\end{array}
\right.
\end{equation*}
Existence and uniqueness of a solution $(K,\Lambda)$ to the backward stochastic  Riccati equation (BSRE)  in \reff{EDSRcommon} is discussed in Section 3.2 in \cite{pha16} by relating  BSRE to standard LQ control problems. Given $(K,\Lambda)$, existence of a solution $(Y,Z^{Y},Z^{Y,0})$ to the linear mean-field BSDE in \reff{EDSRcommon} is obtained as in Step 3(iii) of Section 
\ref{Sec:procedure} by results in  \cite[Thm.~2.1]{LiSunXiong}. 
Finally, the optimal control is given by
\beqs
\alpha_t^* &=& - S_t^{-1} U_t \big(X_t-\E[X_t^*|W^0_t]\big) - S_t^{-1}\big(\xi_t - \E[\xi_t|W^0_t]\big) - \hat S_t^{-1} \big(V_t \E[X_t^*|W^0_t] + O_t \big),
\enqs
where we have set $X^*$ $=$ $X^{\alpha^*}$.
}
\end{remark}

\begin{remark} \label{remcoefsto} 
{\rm 
The method we propose requires some coefficients to be deterministic. Namely, only $\beta, \gamma, M, H, L$ are here allowed to be stochastic. Indeed, if any other coefficient in the SDE or in the cost functional were stochastic, after completing the square we would have a term in the form $\eee[\Xi_t X^\alpha_t]^2$, with $\Xi_t$ stochastic. Due to the randomness of $\Xi_t$, this term cannot be rewritten to match the terms in the candidate $w_t (X^\alpha_t, \eee[X^\alpha_t])$. Conversely, if \textbf{(H1)} holds, $\Xi_t$ is deterministic and the term above rewrites as $\Xi_t \eee[X^\alpha_t]^2$.
\ep
}
\end{remark}

\section{The Infinite-Horizon Problem}  \label{Sec:infinitepb}

We now consider an infinite-horizon version of the problem in \eqref{pb:payoff} and adapt the results to this  framework. The procedure is similar to the finite-horizon case, but non-trivial technical issues emerge when dealing with the equations for $(K,\Lambda,Y,R)$ and the admissibility of the optimal control.

\vspace{0.2cm}


On a filtered probability space $(\Omega,\Fc,\F,\P)$ as in Section \ref{Sec:intro} with $\F$ $=$ $(\Fc_t)_{t\geq 0}$,  let $\rho > 0$ be a positive discount factor, and define the set of admissible  controls as
\beqs
\aaa &:=&  \left\{ \alpha: \Omega \times\R_+ \to \rr^m \text{ s.t.~$\alpha$ is $\F$-adapted and}\int_0^\infty e^{-\rho t} \eee[|\alpha_t|^2] dt < \infty \right\},
\enqs
while the controlled process is defined on $\rr_+$ by
\begin{equation}
\label{pb:SDEINF}
\begin{cases}
dX^\alpha_t = b_t\big(X^\alpha_t, \eee[X^\alpha_t], \alpha_t, \eee[\alpha_t]\big) dt + \sigma_t\big(X^\alpha_t, \eee[X^\alpha_t], \alpha_t, \eee[\alpha_t]\big) dW_t, \quad t \geq 0, 
\\
X^\alpha_0=X_0,
\end{cases}
\end{equation}
where for each $t \geq 0$, $x,\bar x \in \rr^d$ and $a,\bar a \in \rr^m$ we  now set
\begin{equation}
\label{pb:coeffSDEINF}
\begin{aligned}
& b_t\big(x, \bar x, a, \bar a \big) := \beta_t + B x + \tilde B \bar x + C a + \tilde C \bar a,
\\
& \sigma_t\big(x, \bar x, a, \bar a \big) := \gamma_t + D x + \tilde D \bar x + F a + \tilde F \bar a.
\end{aligned}
\end{equation}
Notice that, unlike Section \ref{Sec:intro} and as usual in infinite-horizon problems, the coefficients of the linear terms are constant vectors, but  the coefficients 
$\beta$ and $\gamma$ are allowed to be stochastic processes. 

\vspace{1mm}

The control problem on infinite horizon is formulated as
\begin{equation} \label{pb:payoffINFI}
\begin{array}{ccl}
J(\alpha) &:=&  \eee \Big[ \int_0^\infty e^{-\rho t}  f_t\big(X^\alpha_t, \eee[X^\alpha_t], \alpha_t, \eee[\alpha_t]\big) dt \Big], \\
\rightarrow \;\;\; V_0 &:=&  \Inf_{\alpha \in \aaa}  J(\alpha),
\end{array}
\end{equation}
where, for each $t \geq 0$, $x,\bar x \in \rr^d$ and $a,\bar a \in \rr^m$ we have set
\beq 
f_t\big(x, \bar x, a, \bar a \big) &:=& (x-\bar x)\trans Q (x-\bar x) + \bar x\trans(Q + \tilde Q) \bar x +   2a\trans I (x-\bar x) + 2 \bar a\trans (I+\tilde I) \bar x   \nonumber  \\
& & \;\;\; + \;  (a - \bar a)\trans N (a - \bar a)  + \bar a\trans(N + \tilde N) \bar a +  2M_t\trans x + 2H_t\trans a. \label{pb:coeffPayoffINF}
\enq

Notice that, as usual in  infinite-horizon problems, the coefficients of the quadratic terms are here constant matrices, and the only non-constant coefficients are $H,M$, which may be stochastic processes. Given a normed space $(\M,|.|)$, we set 
\beqs
L^{\infty}(\R_+,\M) &:=& \bigg\{ \phi : \R_+\to \M \text{ s.t.~$\phi$ is measurable and $\textstyle \sup_{t \geq 0} |\phi_t| < \infty$} \bigg\}, \\
L^{2}(\R_+,\M) &:=& \bigg\{ \phi : \R_+ \to \M \text{ s.t.~$\phi$ is measurable and $\textstyle \int_0^\infty e^{-\rho t} |\phi_t|^2 dt  < \infty$} \bigg\}, \\
L^2_{\F}(\Omega \times \R_+,\M) &:=& \bigg\{ \phi : \Omega \times \R_+\to \M \text{ s.t.~$\phi$ is $\F$-adapted and $\int_0^\infty e^{-\rho t} \eee[|\phi_t|^2] dt < \infty$} \bigg\},
\enqs
and ask the following conditions on the coefficients of the problem to hold in the infinite-horizon case.

\begin{itemize}
\item[\textbf{(H1')}] The coefficients in \eqref{pb:coeffSDEINF} satisfy:
	\begin{itemize}
	 \item[(i)] 	$\beta,\gamma \in L^2_{\F}(\Omega \times\R_+,\rr^d)$,
	 \item[(ii)]  $B,\tilde B, D, \tilde D \in \rr^{d\times d}$, $C,\tilde C, F, \tilde F \in  \rr^{d\times m}$.
	 \end{itemize}
 \item[\textbf{(H2')}] The coefficients in \eqref{pb:coeffPayoffINF} satisfy:
	\begin{itemize}
	\item[(i)] $Q, \tilde Q \in  \mathbb{S}^d$,   $N, \tilde N \in  \mathbb{S}^m$, $I, \tilde I \in   \rr^{m\times d}$,
	\item[(ii)] $M \in L^2_{\F}(\Omega \times\R_+,\rr^d)$, $H \in L^2_{\F}(\Omega \times\R_+,\rr^m)$,  
	\item[(iii)]  $ N \; > \;  0$, \, $Q \!-\! I\trans N^{-1} I \; \geq \;  0$,
	\item[(iv)]  $ N \!+\! \tilde N \; > \;  0$, \, $(Q \!+\! \tilde Q) \!-\! (I \!+\! \tilde I)\trans (N \!+\! \tilde N)^{-1} (I \!+\! \tilde I) \; \geq \;  0$.
	\end{itemize}
\item[\textbf{(H3')}]  The coefficients in \eqref{pb:coeffSDEINF} satisfy $\rho >  2 \max \big\{ |B| + |D|^2, \,\, |B+\tilde B| \big\}$.  
\end{itemize}

\vspace{2mm}

Assumptions {\bf (H1')} and {\bf (H2')} are simply a rewriting of {\bf (H1)} and {\bf (H2)} for the case where the coefficients do not depend on the time. A further condition {\bf (H3')}, not present in the finite-horizon case, is here required in order to have a well-defined problem, as justified below.

By {\bf (H1')} and classical results, there exists a unique strong solution $X^\alpha = (X^\alpha_t)_{t \geq 0}$ to the SDE in \eqref{pb:SDEINF}. Moreover, by {\bf (H1')} and 
{\bf (H3')}, standard estimates (see Lemma \ref{lem:est} below) lead to:
\begin{equation}
\label{pb:estimXINF}
\int_0^\infty e^{-\rho t} \eee[|X^\alpha_t|^2] dt \; \leq \; \tilde C_\alpha(1+\eee[|X_0|^2]) \;<\; \infty,
\end{equation}
where $\tilde C_\alpha$ is a constant depending on $\alpha$ $\in$ $\Ac$ only via $\int_0^\infty e^{-\rho t} \eee[|\alpha_t|^2] dt$. Also, by {\bf (H2')} and \eqref{pb:estimXINF}, the problem in \eqref{pb:payoffINFI} is well-defined, in the sense that $J(\alpha)$ is finite for each $\alpha \in \aaa$. 

\begin{Lemma}
\label{lem:est}
Under {\bf (H1')} and {\bf (H3')}, the estimate in \eqref{pb:estimXINF} holds for each square-integrable variable $X_0$ and $\alpha \in  \aaa$. 
\end{Lemma}
\noindent {\bf Proof.} By the It\^o formula and the Young inequality, for each $\eps>0$ we have (using shortened bar notations, see Remark \ref{rem:notations}, e.g., $\bar X$ $=$ $\E[X]$)
\begin{align}
\label{est0}
\frac{d}{dt}e^{-\rho t}|\bar X_t|^2 \; &= \; e^{-\rho t}\big(-\rho |\bar X_t|^2 + 2 \bar b_t \cdot \bar X_t \big) \nonumber\\
&\leq \; e^{-\rho t}\Big[-\rho |\bar X_t|^2 + 2 \Big( |\bar \beta_t| |\bar X_t| + |C+\tilde C| |\bar \h_t||\bar X_t| + \bar X_t'(B+\tilde B) \bar X_t \Big) \Big] \nonumber\\
&\leq \; e^{-\rho t}\Big[ \Big( -\rho + 2|B+\tilde B| + \eps \Big) |\bar X_t|^2 + c_\eps \big( |\bar \beta_t|^2 +|\bar \h_t|^2 \big) \Big],
\end{align}
where $c_\eps>0$ is a suitable constant. We define 
\begin{equation*}
\zeta_\eps \; = \;  \big|\eee[X_0]\big|^2 + c_\eps  \int_0^\infty e^{-\rho t} \eee\big[|\beta_t|^2 + |\alpha_t|^2 \big] dt, 
\qquad\quad
\eta_\eps \; = \;  \rho - 2|B+\tilde B| - \eps, 
\end{equation*}
and notice that $\zeta_\eps$ $<$ $\infty$ by {\bf (H1')} and by $\alpha \in \aaa$, while $\eta_\eps > 0$ for $\eps$ small enough, by {\bf (H3')}. Applying the Gronwall inequality, we then get
\begin{equation}
\label{est1} 
\int_0^\infty e^{-\rho t}\big|\eee[X_t]\big|^2 dt  \; \leq \;  \zeta_\eps \int_0^\infty e^{-\eta_\eps t} dt \; \leq \; c_{\alpha,\eps}(1+\eee[|X_0|^2]),
\end{equation}	
for a suitable constant $c_{\alpha, \eps}$. By similar estimates, we have
\begin{align}
\label{est0bis}
&\frac{d}{dt}\eee\big[e^{-\rho t}|X_t - \bar X_t|^2\big] \nonumber \\
&= \; e^{-\rho t}\eee\big[-\rho |X_t - \bar X_t|^2 + 2( b_t - \bar b_t) \cdot (X_t -  \bar X_t) + |\sigma_t|^2 \big] \nonumber \\
&\leq \; e^{-\rho t}\eee\Big[-\rho |X_t \!-\! \bar X_t|^2 \!+\! 2 \Big( |\beta_t \!-\! \bar \beta_t| |X_t \!-\!  \bar X_t| \!+\! |C| |\alpha_t \!-\! \bar \alpha_t| |X_t \!-\! \bar X_t| \!+\! (X_t \!-\! \bar X_t)' B( X_t \!-\! \bar X_t) \Big) \nonumber \\
&\qquad \; + 2 \Big( |\gamma_t|^2 + |D|^2|X_t - \bar X_t|^2 + |D+\tilde D|^2 |\bar X_t|^2 + |F||\alpha_t|^2 + |\tilde F||\bar \alpha_t|^2 \Big) \Big] \nonumber \\
&\leq \; e^{-\rho t}\eee\Big[ \Big( -\rho + 2|B| + 2|D|^2 + \eps \Big) |X_t - \bar X_t|^2 + \tilde c_\eps \big( |\beta_t|^2 + |\gamma_t|^2 + |\alpha_t|^2 + |\bar X_t|^2 \big) \Big],
\end{align}
where $\tilde c_\eps>0$ is a suitable constant, and hence
\begin{equation}
\label{est2} 
\int_0^\infty e^{-\rho t}\eee[|X_t - \eee[X_t]|^2] dt  \; \leq \;  \tilde c_{\alpha,\eps} (1+\eee[|X_0|^2]),
\end{equation}	
for a suitable $\tilde c_{\alpha,\eps}>0$ (recall that $\int_0^\infty e^{-\rho t}\big|\eee[X_t]\big|^2 dt < \infty$ by \eqref{est1}). The estimate in \eqref{pb:estimXINF} follows by \eqref{est1} and  \eqref{est2}, since
\begin{equation*}
\pushQED{\qed} 
\int_0^\infty e^{-\rho t} \eee[|X_t|^2] dt 
=
\int_0^\infty e^{-\rho t} \eee\big[|X_t - \eee[X_t]|^2\big] dt
-
\int_0^\infty e^{-\rho t} \big| \eee[X_t] \big|^2 dt. \qedhere 
\popQED
\end{equation*}


\vspace{3mm}

\begin{remark}
\label{RemSimplH3}
{\rm 
In some particular cases, the condition in {\bf (H3')} can be weakened. For example, assume that $B \leq 0$ and $B+\tilde B \leq 0$. Then, since $x\trans Bx$, $x\trans (B + \tilde B)x \leq 0$ for each $x$, one term in \eqref{est0} and in \eqref{est1} can be simplified, so that {\bf (H3')} simply writes
\beqs
\rho &>&  2 |D|^2.
\enqs
If, in addition, $\gamma = \tilde D = F=\tilde F=0$, then $|\sigma|^2=|D|^2|X|^2$, so that we do not need the estimate on the square which introduces the factor $2$ in \reff{est0bis}. Correspondingly, {\bf (H3')} simply writes
\begin{equation*}
\pushQED{\qed} 
\rho \; > \;  |D|^2. \qedhere 
\popQED
\end{equation*}
}
\end{remark}

The infinite-horizon version of the verification theorem is an easy adaptation of the arguments in Lemma \ref{prop:verif}.

\begin{Lemma}[\textbf{Infinite-horizon verification theorem}]
\label{prop:verifINF}
Let $\{\Wc_t^\alpha, t \geq 0, \alpha \in \aaa\}$ be a family of $\F$-adapted processes in the form $\Wc_t^\alpha$ $=$ 
$w_t(X^\alpha_t, \eee[X^\alpha_t])$ for some $\F$-adapted random field $\{w_t(x,\bar x), t\geq 0, x,\bar x \in \R^d\}$ satisfying
\beq \label{growthwinfi}
w_t(x,\bar x) & \leq & C( \chi_t + |x|^2 + |\bar x|^2), \;\;\; t \in \R_+, \; x,\bar x \in \R^d, 
\enq
for some positive constant $C$, and non-negative process $\chi$ s.t.~$e^{-\rho t} \E[\chi_t]$ converges to zero as $t \to \infty$,  and such that
\begin{itemize}
\item[(i)] the map $t$ $\in$ $\R_+$ $\longmapsto$ $\E[\Sc_t^\alpha]$, with 
$\Sc_t^\alpha$ $=$ $e^{-\rho t} \Wc_t^\alpha + \int_0^t e^{-\rho s}  f_s\big(X^\alpha_s, \eee[X^\alpha_s], \alpha_s, \eee[\alpha_s]\big) ds$, is non-decreasing for all 
$\alpha \in \aaa$;
\item[(ii)]  the map $t$ $\in$ $\R_+$ $\longmapsto$ $\E[\Sc_t^{\alpha^*}]$ is constant for some $\alpha^* \in \aaa$. 
	\end{itemize}
Then, $\alpha^*$ is an optimal control and $\eee[w_0(X_0,\eee[X_0])]$ is the value of the LQMKV control problem  \eqref{pb:payoffINFI}: 
\begin{equation*}
V_0 \; = \;  \eee[ w_0(X_0,\eee[X_0])] \; = \;  J(\alpha^*).
\end{equation*}
Moreover, any other optimal control satisfies the condition (iii). 	
\end{Lemma}
\noindent {\bf Proof.}
Since the integral in \eqref{pb:estimXINF} is finite,  we have $\lim_{t\rightarrow\infty} e^{-\rho t} \eee[|X_t^\alpha|^2]=0$ for each $\alpha$; by \reff{growthwinfi}, we deduce that 
$\lim_{t\rightarrow\infty} e^{-\rho t} \eee[|\Wc^\alpha_t|]$ $=$  $\lim_{t\rightarrow\infty} e^{-\rho t} \eee[|w_t(X^\alpha_t, \eee[X^\alpha_t])|]$ $=$  $0$. Then, the rest of the proof follows the same arguments as in the one of Lemma \ref{prop:verif}.
\qed

\vspace{2mm}

The steps  to apply Lemma  \ref{prop:verifINF} are the same as the ones in the finite-horizon case. We report the main changes in the arguments with respect to the 
finite-horizon case. 

\vspace{2mm}

\noindent {\bf Steps 1-2.}  
We search for a random field $w_t(x,\bar x)$ in a quadratic form as in \reff{wquadra}: 
\beqs
w_t(x,\bar x) &=& (x - \bar x)\trans K_t (x - \bar x) + \bar x\trans \Lambda_t \bar x + 2 Y_t\trans x + R_t,
\enqs
where the mean optimality principle of Lemma \ref{prop:verifINF} leads now to the following system 
\begin{equation} \label{sysKinfi}
\left\{
\begin{array}{ccl}
dK_t &=&  -  \Phi_t^0(K_t) dt, \;\;\; t \geq 0,   \\ 
d\Lambda_t &=& -  \Psi_t^0(K_t,\Lambda_t)   dt,  \;\;\;  t \geq 0,   \\
dY_t &=& - \Delta_t^0(K_t,\Lambda_t,Y_t,\E[Y_t],Z_t^Y,\E[Z_t^Y]) dt + Z_t^Y dW_t, \;\;  t \geq 0,    \\
dR_t &=& \big[ \rho R_t - \E[\Gamma_t^0(K_t,Y_t, \E[Y_t],Z_t^Y,\E[Z_t^Y]) \big] dt, \;\;\; t \geq 0.  
\end{array}
\right.
\end{equation}
Notice that there are no terminal conditions in the system, since we are considering an infinite-horizon framework. The maps $\Phi^0$, $\Psi^0$, $\Delta^0$, $\Gamma^0$ are defined as in \reff{Phi0}, \reff{PhiK}, \reff{defSUVO}, \reff{defxi}, where the coefficients $B,\tilde B, C,\tilde C, D,\tilde D, F,\tilde F, Q, \tilde Q, N, \tilde N, I, \tilde I$ are now constant. 

\vspace{2mm}

\noindent {\bf Step 3.} We now prove the existence of a solution to the system in \reff{sysKinfi}.

\begin{itemize}
\item[(i)] Consider the ODE for $K$. Notice that the map $k$ $\in$ $\S^d$ $\mapsto$ $\Phi^0(k)$ does not depend on time as the coefficients are constant. We then look for a constant 
non-negative matrix $K$ $\in$ $\S^d$ satisfying $\Phi^0(K)$ $=$ $0$, i.e., solution to
 \beq
 Q -\rho K + KB+ B\trans K + D\trans K D & &  \nonumber \\
-  \;  (I + F\trans K D + C\trans K)\trans(N + F\trans K F)^{-1}(I + F\trans K D + C\trans K) &=& 0.  \label{eqKINF}
 \enq
As in the finite-horizon case, we prove the existence of a solution to \reff{eqKINF} by relating it to a suitable  infinite-horizon linear-quadratic control problem. However, as we could not find a direct result in the literature for such a connection, we proceed through  an approximation argument. For $T \in \rr_+ \cup \{\infty\}$ and $x \in \rr^d$, we consider the following control problems:
\beqs
V^T(x) & : =&  \inf_{\alpha \in \mathcal{A}_T} \eee \Big[ \int_0^T e^{-\rho t} \Big( (\tilde X_t^{x,\alpha})\trans Q \tilde X_t^{x,\alpha} + 2 \alpha_t\trans  I \tilde X_t^{x,\alpha}  
			+ \alpha_t\trans N \alpha_t  \Big) dt \Big],
\enqs
where we have set
$$
\mathcal{A}_T := \Big\{ \alpha: \Omega \times [0,T[ \to \rr^d \text{ s.t.~$\alpha$ is $\F$-progr.~measurable and} \int_0^T e^{-\rho t} \eee[|\alpha_t|^2] dt < \infty \Big\},
$$
and where, for $\alpha$ $\in$ $\Ac_T$, $(\tilde X^{x,\alpha})_{0\leq t\leq T}$ is the solution to 
\beqs
d\tilde X_t &=& (B \tilde X_t + C \alpha_t) dt + (D \tilde X_t + F \alpha_t) dW_t, \qquad \tilde X_0 \; = \; x.
\enqs
The integrability condition $\alpha \in \mathcal{A}_T$ implies that $\int_0^T e^{-\rho t} \eee[|\tilde X^{x,\alpha}_t|^2] dt < \infty$, and so the problems $V_T(x)$ are well-defined for any 
$T \in \rr_+ \cup \{\infty\}$. If $T$ $<$ $\infty$,  as already recalled in the finite-horizon case, {\bf (H1')}-{\bf (H2')} imply that there exists a unique symmetric 
solution $\{K^T_{t}\}_{t \in [0,T]}$ to
\begin{equation} \label{eqKinter}
\left\{
\begin{array}{rcl}
\frac{d}{dt} K_t^T + Q -\rho K_t^T + K_t^TB+ B\trans K_t^T + D\trans K_t^T D & & \\
- (I + F\trans K_t^T D + C\trans K_t^T)\trans(N + F\trans K_t^T F)^{-1}(I + F\trans K_t^T D + C\trans K_t^T) &=& 0, \; t \in [0,T], \\
K_T^T &=& 0, 
\end{array}
\right.
\end{equation}
and that for every $x \in \rr^d$ we have: $V^T(x)$ $=$ $x\trans K_0^T x$. It is easy to check from the definition of $V^T$ that $V^T(x)$ $\rightarrow$ $V^\infty(x)$ as $T$ goes to infinity, from which we deduce that  
\begin{equation*}
V^\infty(x) \; = \;  \lim_{T\rightarrow\infty}  x\trans K_{0}^T x  \; = \;  x\trans(\lim_{T\rightarrow\infty}  K_{0}^T )x, \;\;\; \forall x \in \R^d. 
\end{equation*}
This implies the existence of the limit $K$ $=$ $\lim_{T\rightarrow\infty} K_0^T$.  By passing to the limit in $T$ in the above ODE \reff{eqKinter} at $t$ $=$ $0$, we obtain by standard arguments (see, e.g., Lemma 2.8 in \cite{SunYong}) that  $K$ satisfies \eqref{eqKINF}.  Moreover, $K \in \S^d$ and $K \geq 0$ as it is the limit of symmetric non-negative matrices.
\item[(ii)] Given $K$, we now consider the equation for $\Lambda$.  Notice that the map $\ell$ $\in$ $\S^d$ $\mapsto$ $\Psi^0(K,\ell)$ does not depend on time as the coefficients are constant. We then look for a constant non-negative matrix $\Lambda$ $\in$ $\S^d$ satisfying $\Phi^0(K,\Lambda)$ $=$ $0$, i.e., solution to
\beq
\hat Q^K  -\rho \Lambda  + \Lambda(B+\tilde B) + (B+\tilde B) \trans  \Lambda & & \nonumber \\
- \big(\hat I^K + (C+\tilde C)\trans\Lambda \big)\trans(\hat N^K)^{-1} \big(\hat I^K + (C+\tilde C)\trans \Lambda \big)  &=& 0, \label{eqLINF}
\enq
where we set
\beqs
\hat Q^K &:=& (Q + \tilde Q) + (D+\tilde D)\trans K (D+\tilde D), \\
\hat I^K &:=& (I+\tilde I) + (F + \tilde F)\trans K (D + \tilde D), \\
\hat N^K &:=& (N + \tilde N) + (F + \tilde F)\trans K (F + \tilde F). 
\enqs
Existence of a solution to \eqref{eqLINF} is obtained  by the same arguments used for \eqref{eqKINF} under {\bf (H2')}. 	
\item[(iii)]  Given $(K,\Lambda)$, we consider the equation for $(Y,Z^Y)$. This is a mean-field linear BSDE on infinite-horizon
 \beq 
dY_t &=&  \Big( \vartheta_t + G\trans(Y_t - \eee[Y_t]) + \hat G\trans \eee[Y_t] + J \trans(Z_t^Y - \eee[Z_t^Y]) + \;  \hat{J}\trans\eee[Z_t^Y] \Big) dt  \nonumber \\
& & \hspace{2cm} 	+ \; Z_t^Y dW_t,  \;\;\;\;\; t \geq 0	\label{eqYinfi}
 \enq
where $G$ $\hat G$, $J$, $\hat J$ are constant coefficients in $\R^{d\times d}$, and  $\vartheta$ is a random process in $L^2_{\F}(\Omega\times\R_+,\R^d)$  defined by
\begin{align*}
	& G \; := \;  \rho \,\, \I_d  - B + C S^{-1} U, 
	\\
	& \hat G \; := \;  \rho \,\, \I_d - B - \tilde B + (C +\tilde C)\hat S^{-1} V,
	\\
	& J \; := \;  - D + F S^{-1} U,
	\\
	& \hat{J} \; := \;  - (D + \tilde D) + ( F + \tilde  F) \hat S^{-1} V, 
	\\
	&\vartheta_t \; := \;  - M_t - K (\beta_t - \eee[\beta_t] ) - \Lambda \eee[\beta_t] - D \trans K (\gamma_t - \eee[\gamma_t]) - \hat D \trans K \eee[\gamma_t]
	\\
	& \qquad  \;\;\; + U \trans S^{-1} \big( H_t - \eee[H_t] + F \trans K (\gamma_t - \eee[\gamma_t])\big)  + V \trans \hat S^{-1} \big(\eee[H_t] + (F + \tilde F)\trans K \eee[\gamma_t] \big),
\end{align*}
with
\begin{equation} \label{Sinfi}
\left\{
\begin{array}{rcl}
S  &:=&  N + F \trans K F \; = \; S(K), \\
\hat S &:=& N + \tilde N + (F + \tilde F)\trans K (F + \tilde F)  \; = \; \hat S(K), \\
U &:=& I + F\trans K D + C\trans K \; = \; U(K), \\
V  &:=& I + \tilde I + (F +\tilde F)\trans K (D +\tilde D) +   (C + \tilde C)\trans \Lambda \; = \; V(K,\Lambda). 
\end{array}
\right.
\end{equation}
Although in many practical case an explicit solution is possible (see below), there are no general existence results for such a mean-field BSDE on infinite-horizon, to the best of our knowledge. We then assume what follows.

\begin{itemize}
\item[\textbf{(H4')}] There exists a solution $(Y,Z^Y)$ $\in$  $L^2_{\F}(\Omega\times \R_+,\R^d)\times L^2_{\F}(\Omega\times \R_+,\R^d)$ to \eqref{eqYinfi}. 
\end{itemize}

\vspace{1mm}

\begin{remark}
\label{remH4}
{\rm
In many practical cases, {\bf (H4')}  is satisfied and explicit expressions for $Y$ may be available. We list here  some examples.
\begin{itemize}
\item[-] In the case where $\beta=\gamma=H=M\equiv0$, so that $\vartheta$ $\equiv$ $0$, we see that $Y$  $=$ $Z^Y$ $\equiv 0$ is a solution to \eqref{eqYinfi}, and {\bf (H4')} clearly holds. 
\item[-] If $\beta,\gamma,H,M$ are deterministic (hence, all the coefficients are non-random), the process $Y$ is deterministic as well, that is $Z^Y \equiv 0$ and $\eee[Y]=Y$. Then, the mean-field BSDE becomes a standard linear ODE:
\beqs
dY_t &=& \big( \vartheta_t + \hat G\trans Y_t \big) dt, \;\;\; t \geq 0. 
\enqs
In the one-dimensional case $d=1$, we get
\begin{equation*}
Y_t = - \int_t^\infty e^{- \hat G (s-t)} \vartheta_s ds, \;\;\; t \geq 0.
\end{equation*}
Therefore, when $\hat G - \rho > 0$,  i.e., $(C +\tilde C)\hat S^{-1} V$ $>$ $B+\tilde B$, we have by the Jensen inequality and the Fubini theorem
\begin{equation*}
\int_0^\infty \!\! e^{-\rho t} Y_t^2 dt \leq \tilde c_1 \int_0^\infty \!\! \int_t^\infty \!\! e^{-\rho t}  e^{-\hat G (s-t)} \vartheta_s^2 ds \, dt \leq \tilde c_2 \int_0^\infty \!\! e^{-\rho s} \vartheta_s^2 ds < \infty,
\end{equation*}
for suitable constants $\tilde c_1, \tilde c_2>0$, so that {\bf (H4')} is satisfied. In the multi-dimensional case $d>1$, if $\beta, \gamma,H,M$ are constant vectors (hence, $\vartheta$ is constant as well), we have $Y_t = Y$, with
\begin{equation*}
Y = - (\hat G^{-1}) \trans \vartheta,
\end{equation*}
and {\bf (H4')} is clearly satisfied. 

\item[-] In many relevant cases, the sources of randomness of the state variable and the coefficients in the payoff are independent. More precisely, let us consider a problem with two independent Brownian motions $(W^1,W^2)$ (to adapt the formulas above, we proceed as in Remark \ref{remWmulti}). Assume that only $W^1$ appears in the controlled mean-field  SDE:
\beqs
dX^\alpha_t &=&  \big(\beta_t + B X^\alpha_t + \tilde B \eee[X^\alpha_t] + C \alpha_t + \tilde C \eee[\alpha_t]\big) dt 	\\
& & \;\;\;	+ \; \big(\gamma^1_t + D^1 X^\alpha_t + \tilde D^1 \eee[X^\alpha_t] + F^1 \alpha_t + \tilde F^1 \eee[\alpha_t]\big) dW_t^1,
\enqs
where $\beta,\gamma^1$ are deterministic processes. On the other hand,  the coefficients $M,H$ in the payoff are adapted to the filtration of $W^2$ and independent from $W^1$. In this case, the equation for $\big(Y,Z^Y=(Z^{1,Y},Z^{2,Y}) \big)$ writes
\beqs 
dY_t \!\! &=&  \!\! \Big( \vartheta_t + G\trans(Y_t - \eee[Y_t]) + \hat G\trans \eee[Y_t] + (J^1)\trans(Z_t^{1,Y} - \eee[Z_t^{1,Y}]) + \;  (\hat{J^1})\trans\eee[Z_t^{1,Y}] \Big) dt  \nonumber \\
& & \hspace{2cm} 	+ \; Z_t^{1,Y} dW^1_t + Z_t^{2,Y} dW^2_t,  \;\;\;\;\; t \geq 0,	
 \enqs
where we notice that $Z^{2,Y}$ does not appear in the drift as the corresponding coefficients are zero. Notice that  the process $(\vartheta_t)_t $ is adapted with respect to the filtration of $W^2$, while the other coefficients are constant. Then, it is natural to look for a solution $Y$ which is, as well, adapted to the filtration of $W^2$, i.e., such that $Z^{1,Y}$ $\equiv$ $0$. This leads to the equation: 
\beqs
dY_t &=&  \Big( \vartheta_t + G\trans(Y_t - \eee[Y_t]) + \hat G\trans \eee[Y_t]  \Big) dt  + Z_t^{2,Y} dW^2_t, \;\;\; t \geq 0.
\enqs
For simplicity, let us consider the one-dimensional case $d=1$. Taking expectation in the above equation, we get a linear ODE for $\E[Y_t]$, and a linear BSDE for $Y_t - \eee[Y_t]$, given by
\beqs
d \eee[Y_t]&=&  \big( \E[\vartheta_t] + \hat G \eee[Y_t] \big) dt, \;\;\; t \geq  0,  \\
d(Y_t - \eee[Y_t]) &=& \Big( \vartheta_t - \E[\vartheta_t] + G (Y_t - \eee[Y_t]) \Big) dt +   Z_t^{2,Y} dW^2_t, \;\;\; t \geq 0,
\enqs
which lead to
\begin{equation*}
Y_t  \; = \;  - \int_t^\infty e^{-G(s-t)} \vartheta_s ds - \int_t^\infty \Big(e^{-\hat G(s-t)} - e^{- G(s-t)}\Big) \eee[\vartheta_s] ds.
\end{equation*}
Provided that $G -\rho, \hat G - \rho>0$, and recalling that $\vartheta$  $\in$ $L^2_{\F}(\Omega\times\R_+,\R^d)$, condition {\bf (H4')} is satisfied by the same estimates as above. See \cite{AidBaseiPham} and Section \ref{Sec:applications} for practical examples from mathematical finance with such properties. \qed
\end{itemize}}
\end{remark}

\item[(iv)]  
Given $(K,\Lambda,Y,Z^Y)$ the equation for $R$ is a linear ODE, whose unique solution is explicitly given  by 
\beq \label{eqRinfi}
R_t &=& \int_t^\infty e^{-\rho (s-t)} h_s ds, \;\;\; t \geq 0, 
\enq
where the deterministic function $h$ is defined, for $t$ $\in$ $\R_+$, by
\beqs
h_t &:=&  \eee\big[ - \gamma_t\trans K \gamma_t - \beta_t\trans Y_t - 2\gamma_t\trans Z_t^Y  + \xi_t\trans S^{-1} \xi_t \big] 
- \eee[\xi_t]\trans S^{-1} \eee[\xi_t] + O_t\trans \hat S^{-1} O_t,
\enqs
with 
\begin{equation} \label{xiinfi}
\left\{
\begin{array}{rcl}
\xi_t &:=& H_t + F \trans K \gamma_t  + C \trans Y_t + F\trans Z_t^Y, \\
O_t &:=& \E[H_t] + (F + \tilde F)\trans K \E[ \gamma_t] + (C + \tilde C)\trans\E[Y_t] + ( F+ \tilde F)\trans \E[Z_t^Y].
\end{array}
\right.
\end{equation}
Under assumptions {\bf (H1')}(i), {\bf (H2')}(ii) and {\bf (H4')},  we see that $\int_0^\infty e^{-\rho t} |h_t| dt$ $<$ $\infty$, from which we obtain that 
$R_t$ is well-defined for all $t$ $\geq$ $0$. Therefore,  $e^{-\rho t} |R_t|$ $\leq$ $\int_t^\infty e^{-\rho s} |h_s| ds$ $\rightarrow$ $0$ as $t$ goes to infinity. 
\end{itemize}

\paragraph{Final step.}  
Given $(K,\Lambda,Y,Z^Y,R)$ solution to \reff{eqKINF}, \reff{eqLINF}, \reff{eqYinfi}, \reff{eqRinfi}, i.e.,  to the system in \reff{sysKinfi},  the function 
\beqs
w_t(x,\bar x) &=& (x - \bar x)\trans K (x - \bar x) + \bar x\trans \Lambda \bar x + 2 Y_t\trans x + R_t,
\enqs
satisfies the growth condition \reff{growthwinfi}, and by construction  the conditions (i)-(ii) of the verification theorem in Lemma \ref{prop:verifINF}. 
Let us now  consider as in the finite-horizon case the candidate for the optimal control given by 
\beq
\alpha_t^* &=& a_t^0(X_t^*,\E[X_t^*]) + a_t^1(\E[X_t^*]) \nonumber\\
&=& - S^{-1} U (X_t-\E[X_t^*]) -  S^{-1} \big(\xi_t - \E[\xi_t] \big)  -   \hat S^{-1} \big(V  \E[X_t^*] +  O_t  \big), \;\; t \geq 0,  \label{alphaoptinfi}
\enq
where $X^*$ $=$ $X^{\alpha^*}$ is the state process with the feedback control $a_t^0(X_t^*,\E[X_t^*]) + a_t^1(\E[X_t^*])$, and $S,\hat S,U,V$, $\xi,O$ are recalled in \reff{Sinfi}, \reff{xiinfi}. 
With respect to the finite-horizon case in Proposition \ref{profini}, the main technical issue is to check that $\alpha^*$ satisfies  the admissibility condition in $\Ac$. We need to make an additional condition on the discount factor:

\vspace{5mm}
\noindent \textbf{(H5')} \hspace{0.8cm} $\rho >  2 \max \Big\{ |B-CS^{-1}U| + |D-FS^{-1}U|^2, \,\, |(B+\tilde B)-(C+\tilde C) \hat S^{-1} V| \Big\}$.
\vspace{5mm}

\noindent From the  expression of $\alpha^*$ in \eqref{alphaoptinfi}, we see that $X^*=X^{\alpha^*}$ satisfies, 
\beqs
dX^*_t &:=& b^*_t dt + \sigma^*_t dW_t, \;\;\; t \geq 0, 
\enqs
with
\beqs
b^*_t \; := \;  \beta^*_t + B^* (X^*_t - \eee[X^*_t]) + \tilde B^* \eee[X^*_t], & & 
\sigma^*_t  \; := \;   \gamma^*_t + D^* (X^*_t - \eee[X^*_t]) + \tilde D^* \eee[X^*_t],
\enqs
where we set
\begin{gather*}
B^* \; := \;  B-CS^{-1} U,
\qquad\qquad
\tilde B^* := (B + \tilde B) - (C+\tilde C) \hat S^{-1} V,
\\
D^* \; := \;  D - FS^{-1} U ,
\qquad\qquad
\tilde D^* \; := \;  (D + \tilde D) - (F + \tilde F) \hat S^{-1}V,
\\
\beta^*_t \; := \;  \beta_t  - CS^{-1} (\xi_t-\eee[\xi_t]) - (C + \tilde C) \hat S^{-1} O_t,
\\
\gamma^*_t \; := \;  \gamma_t  - FS^{-1} (\xi_t-\eee[\xi_t]) - (F + \tilde F )\hat S^{-1} O_t.
\end{gather*}
If {\bf (H5')} holds, by replicating the arguments in the proof of Lemma \ref{lem:est} to $X^*$, we get \linebreak $\int_0^\infty e^{-\rho t} \eee[|X^*_t |^2] dt < \infty$, so that $\int_0^\infty e^{-\rho t} \eee[|\alpha^*_t |^2] dt < \infty$ by \eqref{alphaoptinfi}, hence $\alpha^*$ $\in$ $\Ac$. 

\vspace{1mm}

\begin{remark}
\label{RemSimplH5}
{\rm In some specific cases, {\bf (H5')} can be weakened. For example, assume that 
\beqs
B-CS^{-1} U \;\leq \; 0, \qquad  (B + \tilde B) - (C+\tilde C) \hat S^{-1} V \; \leq \;  0, \qquad  \gamma_t = \tilde D = F = \tilde F =0.
\enqs
In this case, the matrices $B^* ,\tilde B^* $ are negative definite and $|\sigma^*_t|^2 = |D| |X^*_t|^2$, so that, by the same arguments as in Remark \ref{RemSimplH3}, Assumption {\bf (H5')} can be simplified into
\begin{equation*}
\pushQED{\qed} 
\rho \;>\; |D|^2. \qedhere 
\popQED
\end{equation*}
}
\end{remark}


To sum up the arguments of this section, we have proved the following result.

\begin{Theorem} \label{thm:optimalINF}
	Under assumptions  {\bf (H1')}-{\bf (H5')}, the optimal control for the LQMKV  pro\-blem on infinite-horizon  \eqref{pb:payoffINFI} is given  by  \reff{alphaoptinfi}, and  
	the corresponding value of the problem is 
	\beqs
	V_0 &=& J(\alpha^*) \; = \;  \eee\big[ (X_0- \eee[X_0])\trans K (X_0- \eee[X_0]) \big] +  \eee[X_0]\trans \Lambda  \eee[X_0] + 2 \eee\big[Y_0\trans X_0] + R_0.
	\enqs
\end{Theorem}

\begin{remark}  \label{remmultiinfi}
{\rm 
The remarks in Section \ref{Sec:remarks} can be immediately adapted to the infinite-horizon framework. In particular, as in Remark \ref{remH2}, 
one can have existence of a solution to \reff{eqKINF}-\reff{eqLINF}, even when condition {\bf (H2')}  is not satisfied, and obtain the optimal control as in \reff{alphaoptinfi} provided that {\bf (H4')}-{\bf (H5')} are satisfied. On the other hand,  the model considered here easily extends to the case where several independent Brownian motions are present, as described in Remark \ref{remWmulti}. Finally, the results can be extended to the common noise case of Remark \ref{remcommonnoise}, recalling that only the coefficients $\beta_t, \gamma_t, M_t, H_t$ are time-dependent and stochastic, namely, adapted to the filtration generated by the pair $(W,W^0)$, where $W^0$ is a Brownian motion independent of $W$. 
}
\qed
\end{remark}

\begin{remark} 
{\rm 
McKean-Vlasov control problems in infinite horizon have also been considered in \cite{HuangLiYong}. Besides the new approach, as outlined in Section \ref{Sec:intro}, the novelty in this paper is the presence of some stochastic coefficients (namely, $\beta_t,\gamma_t,H_t,M_t$). Allowing some coefficients to be stochastic is important from the practical point of view of applications, see the example in the next Section \ref{Sec:applications}.
}
\qed
\end{remark}

\section{Application to Production of Exhaustible Resource}
\label{Sec:applications}


We study an infinite-horizon model of substitutable production goods of exhaustible resource with a large number $N$ of producers,  inspired by  the papers \cite{guelaslio10} and \cite{chasir14}, see also \cite{gra16}.    
For a  producer $i$ $=$ $1,\ldots,N$,  denote by $\alpha_t^i$ her quantity supplied at time $ t \geq 0$, and by $X^i_t$ her current level of reserve in the good.  As in \cite{guelaslio10}, we assume that the dyna\-mics of the reserve is stochastic with a noise proportional to the current level of reserves, hence evolving according to 
\beqs
dX^i_t &=& - \alpha_t^i dt + \sigma X_t^i  dW^i_t,\;\;\;  t \geq 0,  \;\;\; X_0^i \; = \; x_0^i > 0, 
\enqs
where $\sigma$ $>$ $0$, and $W^i$, $i$ $=$ $1,\ldots,N$, are independent standard Brownian motions. The selling price $P^i$ for producer $i$ follows a linear inverse demand rule, as in \cite{chasir14}, and is subject to a permanent market impact depending on the average extracted quantity of the other producers. It is then given by 
\beqs
P^i_t &=& P_t^0 - \delta \alpha_t^i - \eps \int_0^t \frac{1}{N} \sum_{j=1}^N \alpha_s^j ds, 
\enqs
where $\delta, \eps$ $>$ $0$ are positive constants, and $P^0$ is some continuous random process driven by a Brownian motion $W^0$ independent of $W^i$. The interpretation is that the exogenous price $P^0$ in absence of transaction is independent of the idiosyncratic noises of the producers.  We assume that the filtration generated by the  common observation of the process $P^0$ is equal to the natural filtration $\F^0$ of $W^0$. 

The gain functional for producer $i$ is\footnote{We thank Ren\'e Aid for insightful discussions on this example.}
\beqs
J^i(\alpha^1,\ldots,\alpha^N) &:=& \E \Big[ \int_0^\infty e^{-\rho t} \Big(  \alpha_t^i P^i_t  -   \eta {\rm Var}(\alpha_t^i | W^0)  
-  c \alpha_t^i \big( \frac{ x_0^i - X_t^i}{x_0^i} \big)  \Big) dt \Big],
\enqs
where $\rho > 0$ is the discount rate over an infinite horizon.  The first term represents the instantaneous profit from selling quantity $\alpha^i$ at price $P^i$, the second term penalizes via the non-negative parameter  $\eta$ high individual variations of the produced quantity (given the observation of the process $P^0$) measured by the theoretical (conditional) variance, 
while  the last term $\Cc_i(\alpha^i)$ $=$  $c \alpha_t^i\frac{ x_0^i - X_t^i}{x_0^i}$, with 
$c$ $>$ $0$, represents the cost of extraction. In the beginning, this cost is negligible, and increases as the reserve is depleted. Notice that we assume that the constants $c$ and $\eta$ are the same for all the producers $i$, i.e., the producers are indistinguishable.  

We consider a social planner who imposes the same feedback control policy for all the producers $\alpha_t^i$ $=$ ${\bf a}(t,X_t^i,(P_s^0)_{0\leq s\leq t})$ for some 
measurable function ${\bf a}$ on $\R_+\times\R\times C(\R_+;\R)$, and look for a Pareto optimality among all the producers.  This means that, in contrast with Nash equilibrium where  the producers  act strategically, i.e., each control is perturbed one at a time,  here, we focus on a cooperative equilibrium  
where all the controls are perturbed simultaneously. From the theory of propagation of chaos, the individual level of reserve $X^i$ and price process $P^i$, $i$ $=$ $1,\ldots,N$,  
become independent and identically distributed, conditionally on $P^0$, when $N$ goes to infinity, 
with a common distribution given by the law of  the solution $(X,P)$ to the stochastic  McKean-Vlasov equation
\begin{equation} \label{XPreduced}
\left\{
\begin{array}{ccl}
dX_t &=& - \alpha_t dt + \sigma X_t dW_t, \\
P_t  &=& P_t^0 - \delta \alpha_t  - \eps  \int_0^t \E[ \alpha_s | W^0 ] ds, 
\end{array}
\right.
\end{equation}
for some Brownian motion $W$ independent of $W^0$, and where $\alpha_t$ $=$ ${\bf a}(t,X_t,(P_s^0)_{0\leq s\leq t})$, 
$t\geq 0$. We are then reduced to the problem of a representative producer 
with initial reserve $x_0$ $>$ $0$, dynamics of level of resource $X$ as in \reff{XPreduced},  controlled by the extracted quantity $\alpha$, and selling price $P$ as in \reff{XPreduced}. Her objective is to maximize over $\alpha$ $\in$ $\Ac$, i.e., the set of $\R$-valued  progressively measurable process w.r.t.~the natural filtration of $(W,W^0)$, the gain functional 
\beqs
J(\alpha) &:=& \E \Big[ \int_0^\infty e^{-\rho t} \Big(  \alpha_t P_t   -   \eta {\rm Var}(\alpha_t|W^0) - c \alpha_t \big( \frac{ x_0 - X_t}{x_0} \big) \Big) dt \Big]. 
\enqs
By noting that $\E[X_t|W^0]$ $=$ $x_0 - \int_0^t \E[\alpha_s|W^0] ds$,  so that $P_t$ $=$ $P_t^0 - \delta\alpha_t - \eps(x_0 - \E[X_t|W^0])$, we see that we are  in the framework of Section \ref{Sec:infinitepb} with $d=m=1$ (one-dimensional state variable and control) in the common noise case of Remark \ref{remcommonnoise}, and the coefficients in \eqref{pb:coeffSDEINF} and \eqref{pb:coeffPayoffINF} are given by:
\beqs
C = -1, \quad
D =\sigma, \quad
N = \delta + \eta, \quad
N + \tilde N = \delta, \quad \\
I = - \frac{c}{2x_0}, \quad
I + \tilde I = - \frac{c +  \eps x_0}{2x_0}, \quad 
H_t =  \frac{c + \eps x_0- P_t^0}{2},
\enqs
while the other coefficients are identically zero. Notice that $(H_t)_t$ is a random $\F^0$-adapted process.   
Under the following assumptions
\begin{equation} \label{condP0}
\int_0^\infty e^{-\rho t}\E[|P^0_t|^2] dt < \infty, \qquad\qquad \rho > \sigma^2,
\end{equation}
it is clear that {\bf (H1')} and {\bf (H3')} hold true (for the condition in {\bf (H3')}, we can omit the factor $2$, see Remark \ref{RemSimplH3}). The equations for 
$K$ and $\Lambda$ read as
\begin{equation}\label{eqKLex2}
\left\{
\begin{array}{ccc}
 \frac{ (K+ \frac{c}{2x_0})^2}{\delta+\eta}  + (\rho - \sigma^2) K &=& 0, \\
 \frac{ (\Lambda + \frac{c + \eps x_0}{2x_0})^2}{\delta} + \rho \Lambda - \sigma^2 K &=& 0. 
 \end{array}
 \right.
\end{equation}
Notice that condition {\bf (H2')} is not satisfied. However, we have existence of a solution $(K,\Lambda)$  to  \reff{eqKLex2} such that 
$K_\eta$ $:=$ $\frac{K+\frac{c}{2x_0}}{\delta + \eta}$ $>$ $0$, $\Lambda_\eps$ $:=$ $\frac{\Lambda + \frac{c+\eps  x_0}{2x_0}}{\delta}$ $>$ $0$, and given by
\begin{equation}\label{KL}
\left\{
\begin{array}{rcl}
K_\eta  
&=& \frac{- (\rho-\sigma^2) + \sqrt{ (\rho-\sigma^2)^2 + 2c\frac{\rho-\sigma^2}{x_0(\delta+\eta)}}}{2}   \; > \; 0, \\
\Lambda_\eps 
&=&   \frac{- \rho + \sqrt{ \rho^2 + 2\frac{\rho (c+\eps x_0) + 2\sigma^2 K}{\delta x_0}}}{2}  \; > \; 0.
\end{array}
\right.
\end{equation}
The (linear) BSDE  for $Y$ is written as 
\beqs
d Y_t  &=&  \big[ (\rho + \Lambda_\eps) Y_t  - \frac{\Lambda_\eps}{2}( c + \eps x_0 - P_t^0) \big] dt + Z_t^{Y,0} dW_t^0 , \;\;\; t \geq 0, 
\enqs
whose solution is explicitly given by 
\beqs
Y_t &=& - \E \Big[ \int_t^\infty  \Lambda_\eps e^{-(\rho+\Lambda_\eps)(s-t)}  \frac{P_s^0 - c - \eps x_0}{2}   ds \big| \Fc_t^0 \Big], \;\;\; t \geq 0. 
\enqs
It clearly satisfies  condition {\bf (H4')} from the square integrability condition \reff{condP0} on $P^0$.  
We also notice with Remark \ref{RemSimplH5} that the condition in {\bf (H5')} here writes as $\rho > \sigma^2$, which is satisfied. 
By Theorem \ref{thm:optimalINF},  the optimal control is then given by 
\beqs
\alpha_t^* &=& K_\eta(X_t^* - \E[X_t^*|W^0]) + \Lambda_\eps \E[X_t^*|W^0] \\
& &  \;\; + \;   \frac{1}{2\delta} \Big( P_t^0  - \int_t^\infty  \Lambda_\eps e^{-(\rho+\Lambda_\eps)(s-t)} \E[P_s^0|\Fc_t^0] ds   - (c+\eps x_0) \frac{\rho}{\rho+ \Lambda_\eps} \Big) , 
\enqs
with a conditional optimal level of reserve  given  by
\beq
\E[X_t^* | W^0] &=& x_0 e^{-\Lambda_\eps t} +   
\frac{\rho (c +\eps x_0) }{2\delta}  \frac{1- e^{-\Lambda_\eps t}}{\Lambda_\eps(\rho+ \Lambda_\eps)}  \label{meanXex2} \\
& & \; -  \; \frac{1}{2\delta}   \int_0^t e^{-\Lambda_\eps(t-s)} \Big( P_s^0 - \int_s^\infty \Lambda_\eps e^{-(\rho+\Lambda_\eps)(u-s)} \E[P_u^0|\Fc_t^0] du \Big) ds, \;\;\;  t \geq 0. \nonumber
\enq
Suppose that the price $P^0$ admits a stationary level in mean, i.e., $\eee[P_t^0]$ converges to some constant $\bar p$ when $t$ goes to infinity: $\bar p$ is interpreted as a substitute price for the exhaustible resource.  In this case,  it is easy to see from \reff{meanXex2} that the optimal level of reserve also admits a stationary level in mean: 
\beqs \label{explixinfi}
\lim_{t\rightarrow\infty} \E[X_t^*] &=& 
\frac{\rho(c + \eps x_0 - \bar p)}{2\delta\Lambda_\eps(\rho + \Lambda_\eps)} \; =: \; \bar x_\infty.
\enqs
From straightforward algebraic calculations on  \reff{eqKLex2}, we have: 
\beqs
2\delta \Lambda_\eps(\rho + \Lambda_\eps) &=& \rho\eps +  \frac{K_\eta + \rho}{K_\eta + \rho - \sigma^2} (\rho -\sigma^2) \frac{c}{x_0}, 
\enqs 
and thus
\beq \label{xinfty}
\bar x_\infty &=&  \frac{1}{\frac{\eps x_0}{ c + \eps x_0} + \frac{c}{c + \eps x_0} \lambda_\eta}
\big(1 - \frac{\bar p}{c + \eps  x_0}\big) x_0. 
\enq
with $\lambda_\eta$ $:=$ $\frac{\rho-\sigma^2}{\sigma^2}\frac{K_\eta}{K_\eta+\rho-\sigma^2}$ $\in$ $(0,1)$.  
The term $c + \eps x_0$ is the cost of extraction for the last unit of resource. When it is larger than the substitute price $\bar p$, i.e., the Hotelling rent $H_r$ $:=$ 
$\bar p - c - \eps x_0$ is negative, this ensures that the average long term level of reserve $\bar x_\infty$ is positive, 
meaning that there is remaining resource when switching to the substitute good. 

One can study the sensitivity of $\bar x_\infty$ $=$ $\bar x_\infty(\eta,\eps)$ w.r.t.~the intermittence parameter $\eta$ and permanent market impact $\eps$. When the Hotelling rent $H_r$ is negative, and noting that $K_\eta$ is  decreasing with $\eta$, and so $\lambda_\eta$ is increasing with $\eta$,  we see from \reff{xinfty} that $\bar x_\infty$ is decreasing with $\eta$, and 
\beqs
\bar x_\infty & \searrow & x_0 \big(1 - \frac{\bar p}{c + \eps  x_0}\big), \;\;\; \mbox{ as } \; \eta \nearrow \infty. 
\enqs
On the other hand, for fixed $\eta$, we easily see from \reff{xinfty} that
\beqs
\lim_{\eps \rightarrow 0} \bar x_\infty \; = \;  \frac{1}{\lambda_\eta} x_0\big(1 - \frac{\bar p}{c} \big), & \mbox{ and } & \lim_{\eps\rightarrow\infty}  \bar x_\infty \; = \; x_0. 
\enqs
Finally, notice that the existence of a stationary level of resource in mean implies that  $\lim_{t\rightarrow\infty}\E[\alpha_t^*]$ $=$ $0$. In other words, one stops on average to extract the resource in the long term.

\section{Conclusion}
\label{Sec:concl}

In this paper we propose a weak martingale approach to solve linear-quadratic McKean-Vlasov stochastic control problems. In particular, we allow some coefficients to be stochastic. We first consider finite-horizon problems, characterizing the value function and the optimal control through a suitable system of BSDEs and ODEs. Precise conditions are set on the coefficients, ensuring that such a system admits a unique solution. We then extend the results to the case where several Brownian motions and a common noise are present. Infinite-horizon problems are also considered. In this case, additional conditions are required to the coefficients of the problem. We finally provide a detailed example from mathematical finance.

\vspace{0.6cm}
{\small \noindent\textbf{Acknowledgements.} This work is part of the ANR project CAESARS (ANR-15-CE05-0024), and also supported by FiME (Finance for Energy Market Research Centre) and the ``Finance et D\'eveloppement Durable - Approches Quantitatives'' EDF - CACIB Chair.}

\vspace{7mm}


\small

\bibliographystyle{plain}

\begin{thebibliography}{1}
	
\bibitem{AidBaseiPham} Aid, R., Basei, M., Pham, H.: The coordination of centralised and distributed generation. ArXiv:1705.01302 
	
\bibitem{BalataEtAl} Balata A., Hur\'e C., Lauri\`ere M., Pham H., Pimentel I.: A class of finite-dimensional numerically solvable McKean-Vlasov control problems. ArXiv:1803.00445
	
\bibitem{benetal13} Bensoussan,  A., Frehse, J., Yam, P.: Mean Field Games and Mean Field Type Control Theory. Springer Briefs in Mathematics (2013) 
	
\bibitem{benetal11} Bensoussan, A., Sung, K. C. J., Yam, S. C. P., Yung, S. P.: Linear-quadratic mean field games. J. Optim. Theory Appl. \textbf{169}(2), 1--34 (2016)

\bibitem{cardel18} Carmona, R., Delarue, F.: Probabilistic Theory of Mean Field Games with Applications I-II. Springer (2018) 

\bibitem{chasir14} Chan, P., Sircar, R.: Bertrand and Cournot mean field games. Appl. Math. Optim. {\bf 71}(3), 533--569 (2015)

\bibitem{elk81} El Karoui, N.: Les aspects probabilistes du contr\^ole stochastique. 9th Saint Flour Probability Summer School-1979. Lecture Notes in Math. {\bf 876}, Springer, 73--238 (1981)

\bibitem{gra16} Graber, P. J.: Linear Quadratic Mean Field Type Control and Mean Field Games with Common Noise, with Application to Production of an Exhaustible Resource.
Appl. Math. Optim. {\bf 74}(3), 459--486 (2016)
 
\bibitem{guelaslio10} Gu\'eant, O., Lasry, J.-M., Lions, P.-L.: Mean field games and applications. Paris-Princeton Lectures on Mathematical Finance 2010, Springer 205--266 (2011)

\bibitem{HuangLiYong} Huang, J., Li, X., Yong, J.: A Linear-Quadratic Optimal Control Problem for Mean-Field Stochastic Differential Equations in Infinite Horizon. Mathematical Control and Related Fields {\bf 5}(1), 97--139 (2015)

\bibitem{LiSunXiong} Li, X., Sun, J., Xiong, J.: Linear Quadratic Optimal Control Problems for Mean-Field Backward Stochastic Differential Equations. Appl. Math. Optim. (2017)
 
\bibitem{LiSunYong} Li, X., Sun, J., Yong, J.: Mean-Field Stochastic Linear Quadratic Optimal Control Problems: Closed-Loop Solvability. Probab. Uncertain. Quant. Risk, \textbf{1}:2 (2016)
	
\bibitem{pha16} Pham, H.: Linear quadratic optimal control of conditional McKean-Vlasov equation with random coefficients and applications. Probab. Uncertain. Quant. Risk \textbf{1}:7 (2016)

\bibitem{pen92} Peng, S.: Stochastic Hamilton Jacobi Bellman equations.  
SIAM J. Control Optim. {\bf  30}(2), 284--304 (1992)

\bibitem{Sun} Sun, J.: Mean-Field Stochastic Linear Quadratic Optimal Control Problems: Open-Loop Solvabilities. ESAIM Control Optim. Calc. Var. {\bf 23}(3), 1099--1127 (2017)

\bibitem{SunYong} Sun, J., Yong, J.: Stochastic Linear Quadratic Optimal Control Problems in Infinite Horizon. Appl. Math. Optim. \textbf{78}(1), 145--183 (2018)	
	
\bibitem{Yong2013} Yong, J.: A Linear-Quadratic Optimal Control Problem for Mean-Field Stochastic Differential Equations. SIAM J. Control Optim. {\bf 51}(4), 2809--2838 (2013)

\bibitem{YongZhou} Yong, J., Zhou, X.: Stochastic controls, Springer (1999)
\end{thebibliography}

\end{document}